\documentclass[pdflatex,sn-mathphys-num]{sn-jnl}

\usepackage{microtype}

\usepackage{graphicx}
\usepackage{booktabs}
\usepackage{xcolor}

\usepackage{mathtools}
\usepackage{amssymb}
\usepackage{amsthm}
\usepackage{bm}
\usepackage{mleftright}

\usepackage{algorithm}
\usepackage{algpseudocode}
\usepackage[title]{appendix}
\usepackage{placeins}

\usepackage{xurl}
\usepackage{xspace}
\usepackage{cleveref}
\hypersetup{bookmarksdepth=3}

\makeatletter
\providecommand*{\theHALG@line}{}
\renewcommand*{\theHALG@line}{%
  \thealgorithm.\arabic{ALG@line}%
}
\makeatother

\newcommand{\real}{\mathbb{R}}
\DeclareMathOperator{\tr}{tr}
\DeclareMathOperator{\diag}{diag}
\newcommand{\mat}[1]{\boldsymbol{#1}}
\renewcommand{\vec}[1]{\boldsymbol{#1}}
\newcommand{\norm}[1]{\left\lVert #1 \right\rVert}
\newcommand{\Id}{\mathbf{I}}
\DeclareMathOperator{\expect}{\mathbb{E}}
\newcommand{\order}{\mathcal{O}}
\newcommand{\set}[1]{\mathsf{#1}}
\newcommand{\e}{\mathrm{e}}
\renewcommand{\hat}[1]{\widehat{#1}}
\renewcommand{\epsilon}{\varepsilon}
\newcommand{\RPCholesky}{\textnormal{\textsc{RPCholesky}}\xspace}
\newcommand{\krill}{\textnormal{\textsc{KRILL}}\xspace}
\newcommand{\FALKON}{\textnormal{\textsc{FALKON}}\xspace}
\newcommand{\BLESS}{\textnormal{\textsc{BLESS}}\xspace}

\theoremstyle{thmstyleone}
\newtheorem{theorem}{Theorem}

\theoremstyle{thmstyletwo}

\theoremstyle{thmstylethree}
\newtheorem{definition}{Definition}

\begin{document}

\title[Robust, randomized preconditioning]{Robust, randomized preconditioning for kernel ridge regression}

\author[1]{\fnm{Mateo} \sur{D{\'i}az}}\email{mateodd@jhu.edu}

\author[2]{\fnm{Ethan N.} \sur{Epperly}}\email{eepperly@berkeley.edu}

\author[3]{\fnm{Zachary} \sur{Frangella}}\email{zfrangella@alumni.stanford.edu}

\author[4]{\fnm{Joel A.} \sur{Tropp}}\email{jtropp@caltech.edu}

\author*[5]{\fnm{Robert J.} \sur{Webber}}\email{rwebber@ucsd.edu}

\affil[1]{\orgdiv{Department of Applied Mathematics \& Statistics}, \orgname{Johns Hopkins University}, \orgaddress{\city{Baltimore}, \state{MD}, \country{USA}}}

\affil[2]{\orgdiv{Department of Mathematics}, \orgname{University of California Berkeley}, \orgaddress{\city{Berkeley}, \state{CA}, \country{USA}}}

\affil[3]{\orgdiv{Management Science \& Engineering}, \orgname{Stanford University}, \orgaddress{\city{Stanford}, \state{CA}, \country{USA}}}

\affil[4]{\orgdiv{Computing \& Mathematical Sciences}, \orgname{California Institute of Technology}, \orgaddress{\city{Pasadena}, \state{CA}, \country{USA}}}

\affil*[5]{\orgdiv{Department of Mathematics}, \orgname{University of California San Diego}, \orgaddress{\city{La Jolla}, \state{CA}, \country{USA}}}

\abstract{We investigate preconditioned conjugate gradient methods for kernel ridge regression (KRR) problems with a moderate to large number of data points ($10^4 \leq N \leq 10^7$).
We develop and analyze two randomized preconditioners with complementary guarantees.
For full-data KRR, \RPCholesky preconditioning requires $\order(N^2)$ arithmetic operations for fixed accuracy under sufficiently rapid eigenvalue decay of the kernel matrix.
For restricted KRR with $k\ll N$ centers, \krill preconditioning requires $\order((N+k^2)k\log k)$ operations with no eigenvalue-decay assumption.
Experiments on benchmark and scientific data sets demonstrate the robustness of both methods relative to existing preconditioners.
}




\maketitle

\section{Motivation}

Kernel ridge regression (KRR) is widely used in data science and scientific computing, including predicting molecular properties \cite{DBB21,WM21}, but its standard implementation requires solving a dense $N\times N$ linear system.
A direct Cholesky solve requires $
\mathcal{O}(N^3)$ arithmetic operations and $\mathcal{O}(N^2)$ storage, making problems with $ N \geq 10^4$ challenging. 
We therefore seek preconditioners that reduce the iteration count of conjugate gradient (CG) while remaining inexpensive to construct.

We introduce two randomized preconditioners.
For full-data KRR, \RPCholesky
constructs a low-rank approximation of the kernel matrix and, under an eigenvalue-decay condition, yields an $\mathcal{O}(N^2)$ solver (\Cref{thm:a_priori}).
For restricted KRR, \krill uses a sparse sign embedding and provides a condition-number guarantee with no eigenvalue assumption needed (\Cref{thm:a_priori_2}).
\Cref{alg:pcgrpc,alg:krill} summarize the methods, and our experiments show that they are more robust than existing preconditioners across a diverse collection of KRR problems.

\subsection{Notation} \label{sec:notation}

For simplicity, we focus on the real setting, although our work extends to complex-valued kernels without significant modification.
The transpose and Moore--Penrose pseudoinverse of $\mat{A}$ are denoted $\mat{A}^*$ and $\mat{A}^\dagger$.
Double bars $\norm{\cdot}$ indicate the Euclidean norm of a vector or the spectral norm of a matrix.
The matrix condition number is $\kappa(\mat{A}) \coloneqq \lVert \mat{A}\rVert \lVert \mat{A}^{-1} \rVert$.
For a positive-definite matrix $\mat{M}$, the $\mat{M}$-weighted inner product norm is denoted $\lVert \vec{z} \rVert_{\mat{M}} \coloneqq (\vec{z}^*\mat{M}\vec{z})^{1/2}$. The function $\lambda_i(\mat{M})$ outputs the $i$th largest eigenvalue of $\mat{M}$.
The $\ell_1$ norm of the vector $\vec{x}$ is denoted by $\|\vec{x}\|_{\ell_1} = \sum_i |x_i|$. The symbol $|\set{S}|$ denotes the cardinality of a set $\set{S}$.

\section{
Algorithms and
best practices} \label{sec:algorithms}

For the KRR user's convenience, this early section highlights our recommended approaches for solving KRR problems.
Specifically, \Cref{sec:exact_intro} introduces the \RPCholesky preconditioner for solving the full-data KRR problem, while
\Cref{sec:restricted_intro} introduces the \krill preconditioner for solving the restricted  KRR problem.

\subsection{Full-data kernel ridge regression}
\label{sec:exact_intro}

We are given $N$ input--output data pairs $(\vec{x}^{(1)}, y_1), (\vec{x}^{(2)}, y_2), \ldots, (\vec{x}^{(N)}, y_N)$ in $\mathcal{X} \times \mathbb{R}$ for training.
Assume that $K : \mathcal{X} \times \mathcal{X} \to \mathbb{R}$ is a positive-definite kernel function,
and define the $N \times N$ positive-semidefinite kernel matrix $\mat{A}$ with entries $a_{ij} = K(\vec{x}^{(i)}, \vec{x}^{(j)})$.

In full-data kernel ridge regression (KRR), we build a prediction function of the form
\begin{equation*}
    \hat{f}(\vec{x}; \vec{\beta}) = \sum_{i=1}^N \beta_i K(\vec{x}^{(i)}, \vec{x}).
\end{equation*}
The coefficients $\vec{\beta} \in \mathbb{R}^N$ for the function are chosen to minimize the regularized least-squares loss
\begin{equation*}
    L(\vec{\beta}) = \lVert \vec{y} - \mat{A} \vec{\beta} \rVert^2 + \mu \vec{\beta}^\ast \mat{A} \vec{\beta}.
\end{equation*}
Here, $\mu > 0$ is a regularization parameter.
Theoretical literature \cite{caponnetto2006optimal} suggests setting the regularization to $\mu = \mathcal{O}(N^p)$ for a value $p \in (0, 2/3)$ depending on the kernel matrix eigenvalue decay and the smoothness of $f$.
In practice, $\mu$ is typically chosen via cross-validation or a grid search.

Minimizing $L(\vec{\beta})$ is a quadratic optimization problem.
A canonical minimizer is the unique solution of
\begin{equation}
\label{eq:exact_eqn}
    (\mat{A} + \mu \Id) \vec{\beta} = \vec{y}.
\end{equation}
The linear system can be solved in $\order(N^3)$ time using a Cholesky decomposition of $\mat{A} + \mu \Id$ and two triangular solves.
However, when there is a medium or large number of data points ($N \geq 10^4$), we 
instead
recommend solving \eqref{eq:exact_eqn} at a reduced cost using conjugate gradient with \RPCholesky preconditioning, described below.
For a fixed relative accuracy and failure probability, we will prove that the \RPCholesky-preconditioned CG solves the full-data KRR problem in $\order(N^2)$ arithmetic operations,
provided the kernel matrix eigenvalues decay at a sufficiently fast polynomial rate (\Cref{thm:a_priori}).

\subsubsection{\RPCholesky preconditioning}

To build a preconditioner for the full-data KRR equations \eqref{eq:exact_eqn},
we begin with a low-rank approximation $\mat{\hat{A}}$ of the kernel matrix $\mat{A}$.
Then we define the preconditioner
\begin{equation*}
    \mat{P} = \mat{\hat{A}} + \mu \Id
\end{equation*}
Its inverse is applied using the formula in \Cref{alg:pcgrpc}.
There are many ways to obtain a low-rank approximation for this purpose, and we adopt the
\RPCholesky algorithm (\Cref{alg:rpcholesky}, analyzed in \cite{chen2022randomly,epperly2025embracerejectionkernelmatrix}).

\RPCholesky is a variant of the partial Cholesky decomposition that
chooses random pivots at each step according to an evolving probability distribution.
For an input parameter $r \in \mathbb{N}$, the algorithm returns a rank-$r$ approximation in factored form:
\begin{equation*}
    \mat{\hat{A}} = \mat{F} \mat{F}^\ast \quad\text{where}\quad \mat{F} \in \real^{N \times r}.
\end{equation*}
The matrix $\mat{\hat{A}}$ is random because it depends on the choice of random pivots.
The average error of $\mat{\hat{A}}$ is comparable to the best rank-$r$
approximation error, as established in \cite[Thm.~5.1]{chen2022randomly} and \cite[Thm.~4.2]{epperly2025embracerejectionkernelmatrix}.
The algorithm accesses only $(r+1)N$ entries
of the kernel matrix; it uses $\mathcal{O}(rN)$ storage;
and it expends $\mathcal{O}(r^2 N)$ arithmetic operations.

\begin{algorithm}[t]
\caption{{\RPCholesky preconditioning} \label{alg:pcgrpc}}
\begin{algorithmic}[1]
\Require Positive semidefinite matrix $\mat{A} \in \real^{N \times N}$, right-hand-side vector $\vec{y} \in \real^N$, regularization coefficient $\mu$, approximation rank $r$, and relative-residual tolerance $\varepsilon$.
\Ensure Approximate solution $\vec{\beta}_\star$ to $(\mat{A} + \mu \Id) \vec{\beta} = \vec{y}$.
\State $\mat{F} \leftarrow \Call{\RPCholesky}{\mat{A},r,\min(100,\lceil r/10 \rceil)}$.
\Comment{See \Cref{alg:rpcholesky}.}
\State $(\mat{U}, \mat{\Sigma}, \sim) \leftarrow \textsc{EconomySizeSVD}(\mat{F})$. \Comment{$\mat{\hat{A}} = \mat{U} \mat{\Sigma}^2 \mat{U}^*$.}
\State Define primitives for preconditioned conjugate gradient:
\begin{align*}
    \texttt{Product} &: \vec{\beta} \mapsto \mat{A}\vec{\beta}  + \mu \vec{\beta}. \\
    \texttt{Preconditioner} &: \vec{\beta} \mapsto \mat{U} \left[ (\mat{\Sigma}^2 + \mu \Id)^{-1}-\mu^{-1}\Id \right]\mat{U}^* \vec{\beta} + \mu^{-1}\vec{\beta}.
\end{align*}
\State $\vec{\beta}_\star \leftarrow \Call{PCG}{\texttt{Product}, \vec{y}, \varepsilon, \texttt{Preconditioner}}$.
\end{algorithmic}
\end{algorithm}

We consider approximation ranks of the form $r=c\sqrt{N}$, where $c$ is independent of $N$.
Constructing and factorizing the resulting preconditioner requires $\mathcal{O}(N^2)$ arithmetic operations.
If such a rank suffices to control the conditioning, then preconditioned CG (PCG) reaches relative energy-norm error $\varepsilon$ in $\mathcal{O}(\log(1/\varepsilon))$ iterations.
Each iteration requires $\mathcal{O}(N^2)$ operations due to multiplication by the dense kernel matrix.
Thus, the total cost is $\mathcal{O}\bigl(N^2\log(1/\varepsilon)\bigr)$ or $\mathcal{O}(N^2)$ for a fixed accuracy $\varepsilon$.
The pseudocode for \RPCholesky preconditioning is given in \Cref{alg:pcgrpc}.

\subsubsection{Empirical performance}

We compare the empirical performance of \RPCholesky and four other preconditioners across 20 regression and classification problems described in \Cref{tab:datasets} of \Cref{sec:data}. 
For each problem, we randomly subsample $N = 1.5 \times 10^4$ data points.
We standardize the features (subtract the mean, divide by the sample standard deviation)
and measure similarity using the squared exponential kernel
\begin{equation}
\label{eq:squared_exponential}
    K(\vec{x}, \vec{y}) = \exp\bigl(-\tfrac{1}{2\sigma^2} \lVert \vec{x} - \vec{y} \rVert^2 \bigr)
    \quad\text{with bandwidth $\sigma = \sqrt{d}$.}
\end{equation}
The bandwidth $\sigma=\sqrt d$ accounts for dimensionality and is on the same scale as the median-distance bandwidth used in many Gaussian KRR applications \cite{huang2024quasi,shi2023kernelanchor}.
This bandwidth scaling can make it challenging to form an effective preconditioner \cite{zhao2024adaptive}.
By contrast, the kernel matrix becomes easier to approximate in the extreme bandwidth regimes: it approaches the identity as $\sigma\to0$ and the all-ones matrix as $\sigma\to\infty$.

\begin{figure}[t]
    \centering
    \begin{minipage}[b]{0.48\textwidth}
        \centering
        \quad
        Larger regularization: $\mu = 10^{-7} N$
        
        \vspace{0.5em}
        \includegraphics[width = \textwidth]{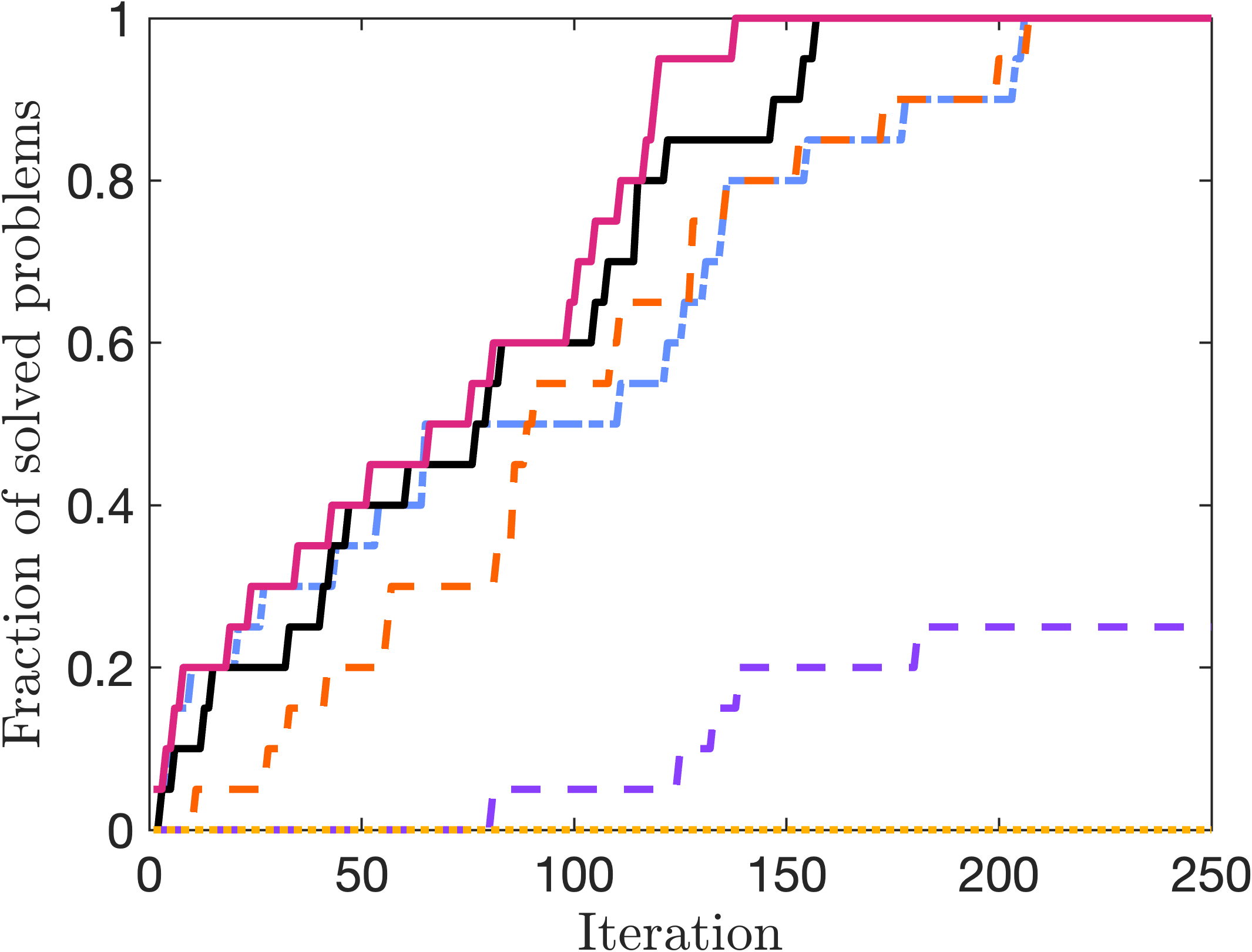}
    \end{minipage}
    \begin{minipage}[b]{0.48\textwidth}
        \centering
        Smaller regularization: $\mu = 10^{-8} N$
        
        \vspace{0.5em}
        \includegraphics[width = \textwidth]{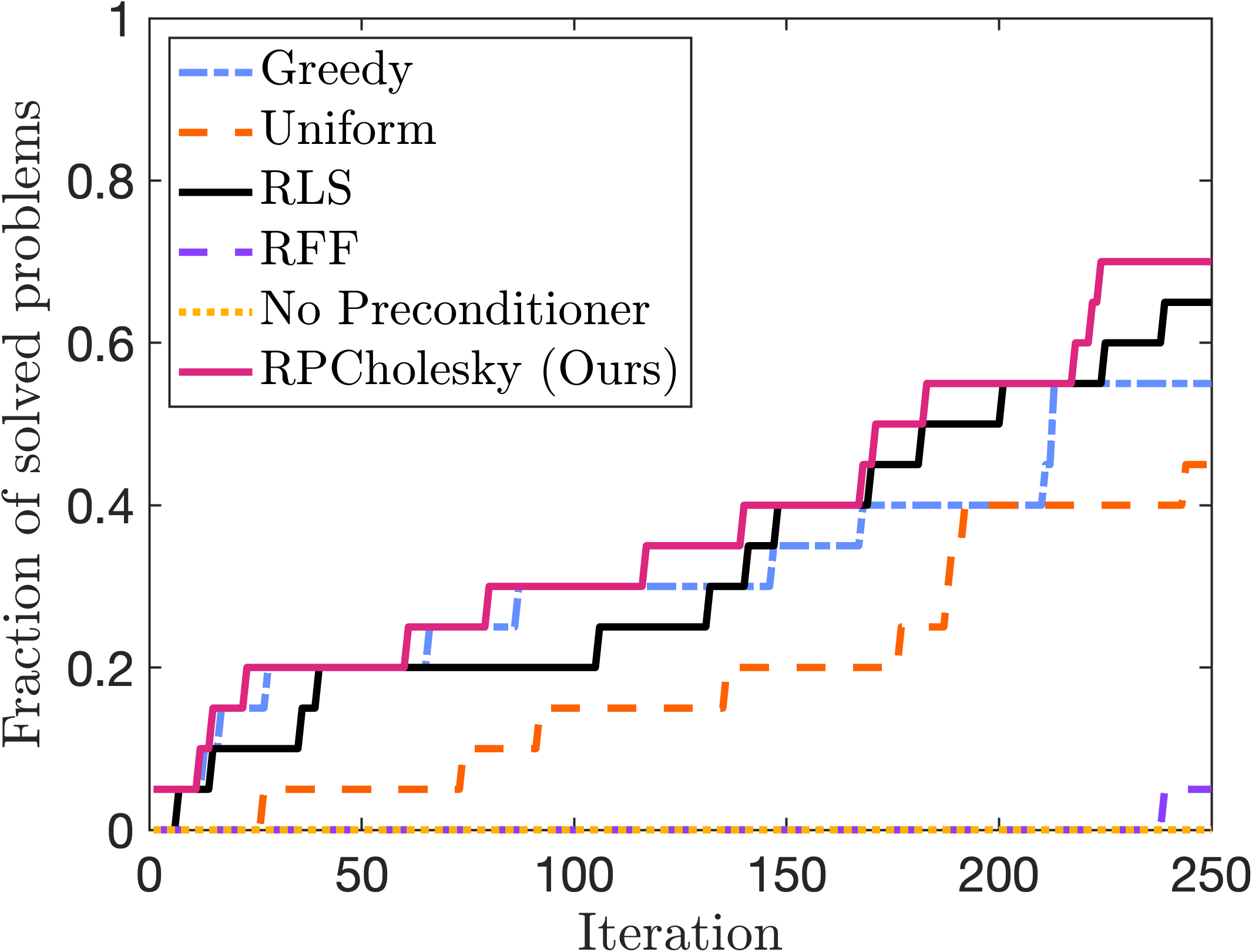}
    \end{minipage}
    \caption{Fraction of solved problems versus number of CG iterations for the 20 KRR problem instances in \Cref{tab:datasets}.}
    \label{fig:exact-performance}
\end{figure}

\begin{figure}[t]
    \centering
    \begin{minipage}[b]{0.48\textwidth}
        \centering        
        \quad \texttt{COMET\_MC\_SAMPLE}
        \includegraphics[width = \textwidth]{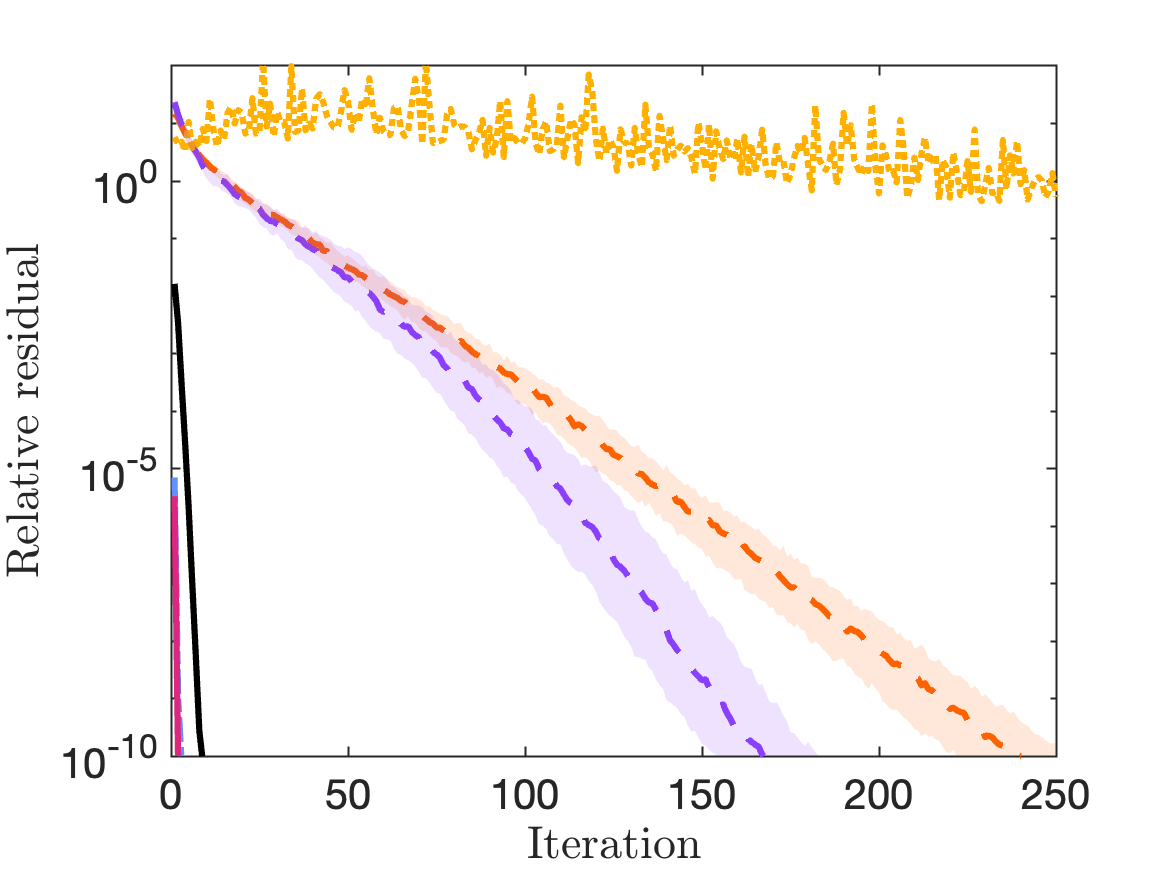}
    \end{minipage}
    \begin{minipage}[b]{0.48\textwidth}
        \centering
        \quad \texttt{w8a}
        \includegraphics[width = \textwidth]{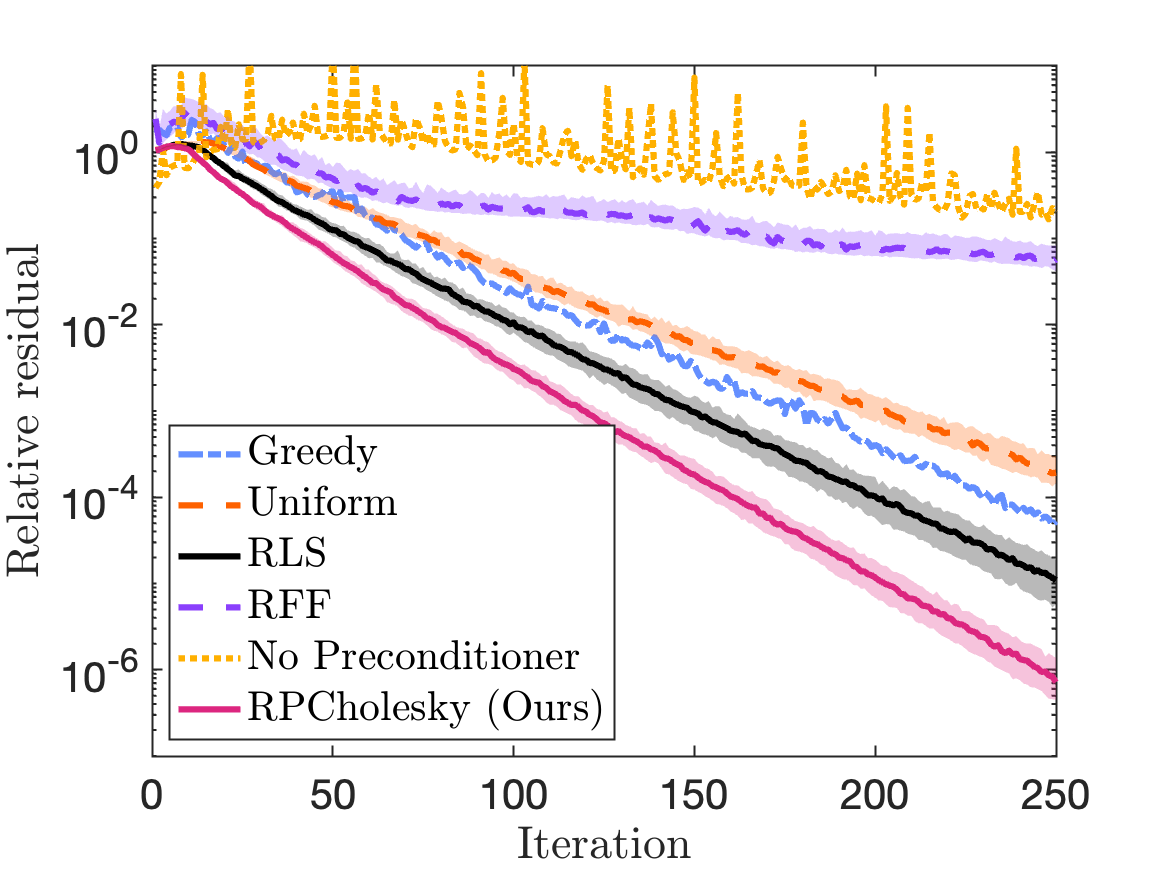}
    \end{minipage}
    \caption{Relative residual versus number of CG iterations for the problem with the fastest (\texttt{COMET\_MC\_SAMPLE}, left)
    and slowest (\texttt{w8a}, right) CG convergence when $\mu/N = 10^{-7}$.
    Note the different vertical axis scales.}
    \label{fig:slower-methods}
\end{figure}

We formulate each KRR problem~\eqref{eq:exact_eqn} with a larger 
regularization parameter 
$\mu/N = 10^{-7}$ (left panel)
or a smaller regularization parameter
$\mu/N = 10^{-8}$ (right panel).
For each problem, we choose the approximation rank $r=1000 \approx 8.2\sqrt{N}$ and run \RPCholesky with the block size $B=100$, matching \Cref{alg:pcgrpc} and the released implementation.
This choice of parameters is typical for full-data KRR problems \cite{avron2017faster,FTU21}, and the parameters are applied systematically over the 20 regression and classification problems without additional tuning.
The smaller value of $\mu$ makes many of these problems highly ill-conditioned.
We run 250 CG iterations, stopping at the first iterate $\vec{\beta}^{(t)}$ for which $\lVert (\mat{A} + \mu \Id) \vec{\beta}^{(t)} - \vec{y} \rVert / \lVert \vec{y} \rVert \leq \varepsilon = 10^{-3}$.
In \Cref{fig:exact-performance,fig:slower-methods}, an instance is called ``solved'' if this criterion is satisfied within 250 iterations.
The choice of relative tolerance $\varepsilon$ is justified by the fact that the test error plateaus after reaching this tolerance for all the problems we considered.

\Cref{fig:exact-performance} charts the fraction of solved problems using \RPCholesky and four other preconditioners based on low-rank approximation.
Two alternative preconditioners use the partial Cholesky decomposition with greedy pivot selection \cite{gardner2018gpytorch,wang2019exact} or uniformly random pivot selection \cite{COCF16,FTU21}.
The final two preconditioners generate randomized low-rank approximations with ridge leverage score sampling (RLS, \cite{alaoui2015fast,musco2017recursive,RCCR19}) or random Fourier features (RFF, \cite{COCF16,avron2017faster}).
For details of these methods, see \Cref{sec:exact_krr_review}.

\RPCholesky is the top-performing method in \Cref{fig:exact-performance} with both the larger (left) and the smaller (right) regularization parameters.
The closest competitor is RLS, but \RPCholesky consistently performs as well or better than RLS.
When $\mu/N = 10^{-7}$, \RPCholesky converges more quickly than RLS on 19 of the 20 test problems and it achieves a speedup greater than $>1.5\times$ for three problems.
As another advantage, \RPCholesky is simpler to implement and faster to run than RLS, as it eliminates the cost of approximating the ridge leverage scores (see \Cref{sec:exact_krr_review}).

To assess sensitivity to the random seed, we repeat each experiment 100 times.
\Cref{fig:slower-methods} shows the distribution of the relative residual $\lVert (\mat{A} + \mu \Id) \vec{\beta} - \vec{y} \rVert / \lVert \vec{y} \rVert$ for the two problems with the fastest (\texttt{COMET\_MC\_SAMPLE}) and the slowest (\texttt{w8a}) CG convergence.
Lines indicate the median error and shaded regions indicate the $20\%$ to $80\%$ error quantiles, which are computed over 100 independent runs.
The error quantiles are hardly distinguishable from the median lines, 
and they are even closer to the median for \RPCholesky than the other stochastic methods (uniform, RLS, and RFF).

The importance of eigenvalue decay is also evident in
\Cref{fig:slower-methods}.
In the problem with rapid eigenvalue decay (\texttt{COMET\_MC\_SAMPLE}, left), \RPCholesky converges to the accuracy threshold $\varepsilon = 10^{-3}$ in just a single iteration.
In the problem with slow eigenvalue decay (\texttt{w8a}, right),
\RPCholesky requires approximately $120$ times as many iterations.
All the methods based on low-rank approximation struggle to control the condition number for the \texttt{w8a} problem.
Nonetheless, even for this difficult problem, \RPCholesky significantly improves on unpreconditioned CG.

\subsubsection{Theoretical guarantees}

\RPCholesky is guaranteed to solve any full-data KRR problem given sufficiently fast eigenvalue decay in the kernel matrix $\mat{A}$.
To explain the eigenvalue decay condition in detail, we introduce a quantitative measure called the \emph{$\mu$-tail rank}:
\begin{definition}[Tail rank]
\label{def:d_tail}
The $\mu$-tail rank of a psd matrix $\mat{A} \in \mathbb{R}^{N \times N}$ is
\begin{equation*}
    \textup{rank}_\mu(\mat{A}) := \min\Biggl\{r \geq 0:\,
    \sum_{i > r} \lambda_i(\mat{A}) \leq \mu\Biggr\}.
\end{equation*}
\end{definition}
Next, we state our main performance guarantee for \RPCholesky preconditioning;
the proof appears in \Cref{sec:a_priori}.
\begin{theorem}[\RPCholesky preconditioning]
\label{thm:a_priori}
Fix a failure probability $\delta \in (0,1)$ and an error tolerance $\epsilon \in (0,1)$.  Let $\mat{A}$ be any positive-semidefinite matrix, and let $\mu > 0$ be a positive number.
Construct a random approximation $\mat{\hat{A}}$ using \RPCholesky with block size $B = 1$ and approximation rank
\begin{equation*}
    r \geq \textup{rank}_\mu(\mat{A}) \,\bigl(1 + \log \bigl(\tr \mat{A} / \mu \bigr) \bigr).
\end{equation*}
With probability at least $1 - \delta$,
the \RPCholesky preconditioner $\mat{P} = \mat{\hat{A}} + \mu \Id$
controls the condition number at a level
\begin{equation} \label{eqn:rpc-success-event}
    \kappa\bigl(\mat{P}^{-1/2}(\mat{A} + \mu \Id) \mat{P}^{-1/2}\bigr) \leq 3 / \delta.
\end{equation}
Conditional on the event~\eqref{eqn:rpc-success-event},
when we apply preconditioned CG to the KRR linear system $(\mat{A} + \mu \Id) \vec{\beta} = \vec{y}$,
we obtain an approximation $\vec{\beta}^{(t)}$ to the actual solution $\vec{\beta}$ 
that satisfies
\begin{equation*}
    \lVert \vec{\beta}^{(t)} - \vec{\beta} \rVert_{\mat{A} + \mu \Id}
    \leq \epsilon\, \lVert \vec{\beta} \rVert_{\mat{A} + \mu \Id}
\end{equation*} 
at any iteration $t \geq \delta^{-1 \slash 2} \log (2 \slash \epsilon)$.
\end{theorem}

\Cref{thm:a_priori} provides a sufficient condition for the approximation rank, which is derived from the trace-error guarantee for \RPCholesky \cite{chen2022randomly}.
The term $\operatorname{rank}_{\mu}(\mat A)$ naturally appears in the rank requirement, since eigenvalues beyond this rank sum to at most $\mu$.
The additional term $1+\log(\tr(\mat{A})/\mu)$ is an oversampling factor from the general \RPCholesky analysis, and we do not know whether it is necessary for preconditioning. 
For comparison, an ideal preconditioner formed from the exact leading eigenvectors would achieve a fixed condition-number bound once $\lambda_{r+1}(\mat{A})=\mathcal{O}(\mu)$, so the optimal rank may be smaller than the sufficient bound above.

\Cref{thm:a_priori} assumes \RPCholesky is performed with a block size $B = 1$ for simplicity, but the theorem can be extended to $B > 1$ given a slightly higher approximation rank: 
\begin{equation*}
r \geq \bigl(\textup{rank}_\mu(\mat{A}) + B\bigr) \,\bigl(1 + \log \bigl(\tr \mat{A} / \mu \bigr) \bigr).
\end{equation*}
The recent work \cite{epperly2025embracerejectionkernelmatrix} analyzes \RPCholesky with arbitrary block size $B \geq 1$, and it introduces an ``accelerated \RPCholesky'' algorithm that improves the accuracy of block \RPCholesky.

If $\operatorname{rank}_\mu(\mat A) \bigl(1+\log(\tr(\mat A)/\mu)\bigr)=\order(\sqrt N)$,
then \Cref{thm:a_priori} guarantees an $\order(N^2)$ solution cost. Uniform and greedy Nystr\"om preconditioning admit no similar guarantee (see \Cref{sec:exact_krr_review}).
\Cref{thm:a_priori} also provides insight into the case where the $\mu$-tail rank is larger than $\order(\sqrt{N})$.  In this case, we may need to use a larger approximation rank $r$, and the construction cost for the preconditioner would exceed $\order(N^2)$ operations.

\begin{figure}[t]
    \centering
    \qquad Original eigenvalues
    \qquad \RPCholesky eigenvalues
    \qquad Shifted eigenvalues

    \qquad \quad $\lambda_i(\mat{A} + \mu \Id)$
    \quad \qquad $\lambda_i(\mu \mat{P}^{-1/2} (\mat{A} + \mu \Id) \mat{P}^{-1/2})$
    \quad \qquad $\lambda_{i+r}(\mat{A} + \mu \Id)$

    \vspace{0.5em}
    \includegraphics[width = \textwidth]{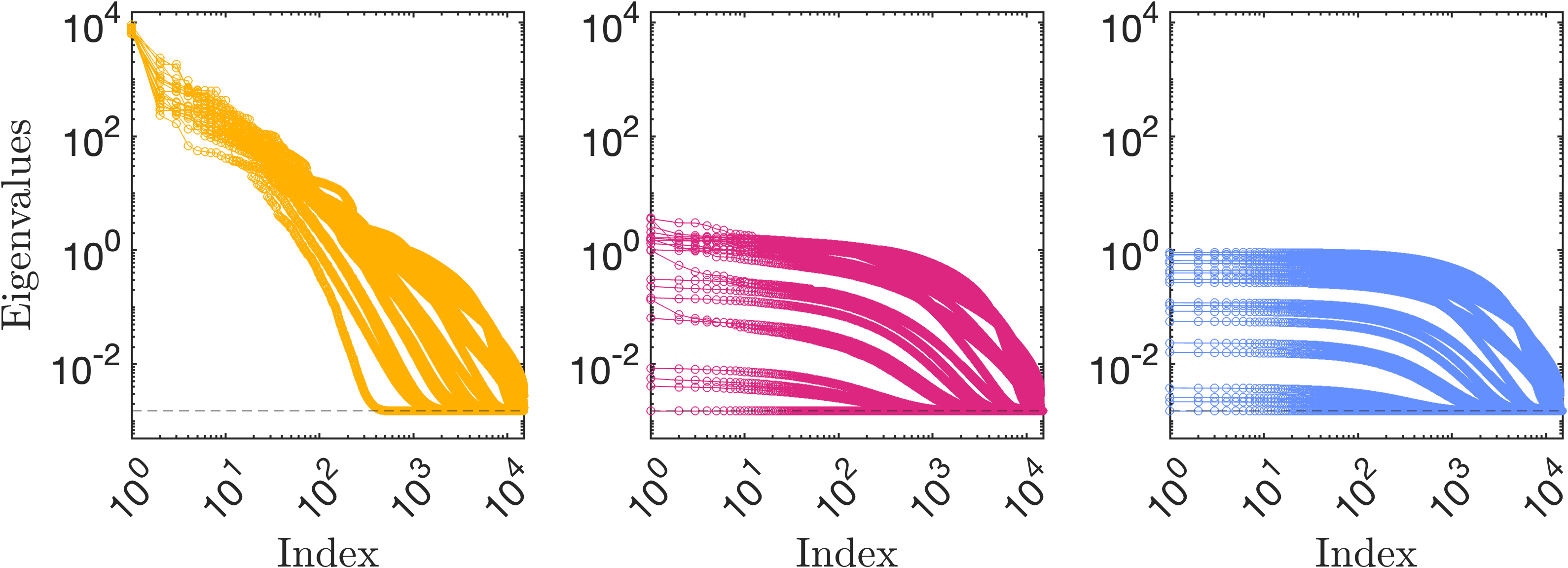}
    \caption{Eigenvalue decay for the 20 KRR problems in \Cref{tab:datasets}.
    Left panel shows eigenvalues of the matrix $\mat{A} + \mu \Id$ before \RPCholesky preconditioning.
    Middle panel shows eigenvalues of the matrix $\mu \mat{P}^{-1/2} (\mat{A} + \mu \Id) \mat{P}^{-1/2}$ after \RPCholesky preconditioning with $r = 1000$.
    Right panel shows eigenvalues $\lambda_{i+r}(\mat{A} + \mu \Id)$ resulting from the mathematically optimal rank-$r$ preconditioner.
    Dashed lines indicate the regularization parameter $\mu = 10^{-7} N$.}
    \label{fig:eigenvalues}
\end{figure}

\Cref{fig:eigenvalues} evaluates \RPCholesky's impact on the 20 example problems when $\mu/N = 10^{-7}$.
Before \RPCholesky preconditioning, the condition number $\lambda_{\max}/\lambda_{\min}$ is nearly $10^7$ (left panel).
\RPCholesky reduces the condition number to a level $10^0$--$10^4$ (middle panel).
The impact of \RPCholesky is similar to the optimal rank-$r$ preconditioner, which eliminates the leading $r$ eigenvalues by setting them equal to $\mu$ (right panel).
Yet even the optimal rank-$r$ preconditioner is sensitive to slow kernel matrix eigenvalue decay.

The optimal use of \RPCholesky balances the cost of forming the preconditioner with the cost of the resulting PCG iterations.
As a simple default, we recommend running \RPCholesky with $r = 10 \sqrt{N}$.
This default is $1.2\times$ larger than the approximation rank used in \Cref{fig:eigenvalues} and is large enough to ensure that all 20 test problems are solved in fewer than 120 iterations when $\mu/N = 10^{-7}$.
See \Cref{sec:homo_energy} for more exploration of the parameter $r$ with a scientific data set.

\subsection{Restricted kernel ridge regression} \label{sec:restricted_intro}

If the number of data points is so large that we cannot apply full-data KRR,
we pursue an alternative approach that we call ``restricted KRR'', which was proposed in \cite{smola2000sparse}.
In restricted KRR, we build a prediction function
\begin{equation*}
    \hat{f}(\vec{x}; \vec{\hat{\beta}}) = \sum_{i=1}^k \hat{\beta}_i K(\vec{x}^{(s_i)}, \vec{x}),
\end{equation*}
using a subset $\vec{x}^{(s_1)}, \vec{x}^{(s_2)}, \ldots, \vec{x}^{(s_k)}$ of input points, which are called ``centers''.  There are many strategies for selecting the centers, such as uniform sampling, ridge leverage score sampling, and \RPCholesky.  One must balance the computational cost of the center selection procedure against the quality of the centers.
For all the experiments in this paper, we use the computationally trivial approach of sampling centers uniformly at random.
However, the development of fast procedures for identifying high-quality centers is a topic for future work.

In restricted KRR, the coefficients $\vec{\hat{\beta}}$ are chosen to minimize the loss function
\begin{equation*}
    L(\vec{\hat{\beta}}) = \lVert \vec{y} - \mat{A}(:, \set{S}) \vec{\hat{\beta}} \rVert^2 + \mu \vec{\hat{\beta}}^\ast \mat{A}(\set{S}, \set{S}) \vec{\hat{\beta}},
\end{equation*}
where $\mat{A}(:, \set{S})$ denotes the submatrix of $\mat{A}$ with column indices in $\set{S} = \{s_1, \ldots, s_k\}$
and $\mat{A}(\set{S}, \set{S})$ denotes the submatrix with row and column indices both in $\set{S}$.
Minimizing this quadratic loss function leads to a coefficient vector $\vec{\hat{\beta}}$ that satisfies
\begin{equation}
\label{eq:restricted_krr_eqn}
    \bigl[\mat{A}(\set{S}, :) \mat{A}(:, \set{S}) + \mu \mat{A}(\set{S}, \set{S})\bigr]\, \vec{\hat{\beta}} = \mat{A}(\set{S}, :) \vec{y}.
\end{equation}
Equation~\eqref{eq:restricted_krr_eqn} is a $k \times k$ linear system that is potentially much smaller than the $N \times N$ system \eqref{eq:exact_eqn} in full-data KRR.
This system can be arbitrarily ill-conditioned and may even be singular.
Therefore, to accurately solve restricted KRR problems in finite-precision arithmetic, we add a small multiple of the machine precision to the diagonal entries of the regularizer $\mu \mat{A}(\set{S}, \set{S})$.

Solving \eqref{eq:restricted_krr_eqn} by direct methods is expensive
because of the $\order(k^2N)$ cost of forming the Gram matrix $\mat{G} = \mat{A}(\set{S}, :) \mat{A}(:, \set{S})$.
When the number of centers is $k \geq 10^3$, we recommend a cheaper approach for solving \eqref{eq:restricted_krr_eqn} using conjugate gradient with the \krill preconditioner, described below.
We will prove that \krill preconditioning solves every restricted KRR problem to any fixed relative accuracy in $\order ((N + k^2)\, k \log k)$ arithmetic operations (\Cref{thm:a_priori_2}).

\subsubsection{\krill preconditioning}
\label{sec:krill_construction}

\begin{algorithm}[t]
\caption{\krill preconditioning \label{alg:krill}}
\begin{algorithmic}[1]
\Require Positive semidefinite matrix $\mat{A} \in \real^{N\times N}$, right-hand-side vector $\vec{y}\in\real^N$, regularization coefficient $\mu$,
centers $\set{S}$, and tolerance $\epsilon$.
\Ensure Approximate solution $\hat{\vec{\beta}}_{\ast}$ to $( \mat{A}(\set{S},:)\mat{A}(:,\set{S}) + \mu \mat{A}(\set{S},\set{S}))\vec{\hat{\beta}} =  \mat{A}(\set{S},:)\vec{y}$.
\State $\mat{H} \leftarrow \mu \mat{A}(\set{S}, \set{S})$
\State $\mat{H} \leftarrow \mat{H} + N  \varepsilon_{\rm mach} \tr(\mat{A}(\set{S}, \set{S}))\, \Id$.
\Comment{$\varepsilon_{\rm mach}\approx 2\times 10^{-16}$ in double precision.}
\State $k \leftarrow |\set{S}|$.
\State $\mat{\Phi} \leftarrow \Call{SparseSignEmbedding}{d=2k,N,\zeta=\lceil\log(k+1)\rceil}$.
\State $\mat{B} \gets \mat{\Phi} \mat{A}(:,\set{S})$.
\State $\mat{P} \leftarrow \mat{B}^\ast \mat{B} + \mat{H}$.
\State $\mat{C} \leftarrow \Call{Cholesky}{\mat{P}}$.
\State Define primitives for preconditioned conjugate gradient:
\begin{align*}
    \texttt{Product} &: \hat{\vec{\beta}} \mapsto
    \mat{A}(\set{S},:)
    \bigl(\mat{A}(:,\set{S})\hat{\vec{\beta}}\bigr)
    +\mat{H}\hat{\vec{\beta}}. \\
    \texttt{Preconditioner} &: \hat{\vec{\beta}} \mapsto \mat{C}^{-1} (\mat{C}^{-\ast} \hat{\vec{\beta}}).
\end{align*}
\State $\vec{\hat{\beta}}_\star \leftarrow \Call{PCG}{\texttt{Product}, \mat{A}(\set{S}, :) \vec{y}, \varepsilon, \texttt{Preconditioner}}$.
\end{algorithmic}
\end{algorithm}

\krill is based on a randomized approximation of the Gram matrix $\mat{G} = \mat{A}(\set{S}, :) \mat{A}(:, \set{S})$.
To form this approximation, we generate a sparse sign embedding $\Phi\in\mathbb{R}^{d\times N}$. 
Independently for each column, we sample $\zeta$ row indices uniformly without replacement, and we place independent Rademacher variables scaled by $\zeta^{-1/2}$ in the corresponding entries; all other entries are zero \cite[\S 9.2]{MT20a}.
Next, we construct the matrix $\mat{B} = \mat{\Phi} \mat{A}(:, \set{S})$.
Finally, to form the \krill preconditioner, we calculate the approximate Gram matrix $\mat{\hat{G}} = \mat{B}^* \mat{B}$ and add the regularizer $\mu\mat{A}(\set{S}, \set{S})$:
\begin{equation*}
    \mat{P}= \mat{B}^* \mat{B} + \mu \mat{A}(\set{S}, \set{S}) = (\mat{\Phi} \mat{A}(:,\set{S}))^* (\mat{\Phi} \mat{A}(:,\set{S})) + \mu \mat{A}(\set{S}, \set{S}).
\end{equation*}
The pseudocode for \krill preconditioning is provided in \Cref{alg:krill}.

\krill preconditioning is efficient
because evaluating the approximate Gram matrix $\mat{\hat{G}}$ requires fewer operations than evaluating the exact Gram matrix $\mat{G}$.
The precise number of operations depends on the sparsity $\zeta$ and the embedding dimension $d$.
In our theoretical analysis, we assume $\zeta$ and $d$ are functions given in \cite{Coh16} that grow according to $\zeta = \order(\log k)$ and $d = \order(k \log k)$.
In practice, we set $\zeta = \lceil\log(k + 1)\rceil$ and $d = 2k$, which we have found to be sufficient to obtain fast CG convergence on all instances.
Recent research \cite{tropp2025comparisontheoremsminimumeigenvalue,chenakkod2025optimalsubspaceembeddingsresolving} has improved the theoretical understanding of sparse random sign embeddings, but the optimal sparsity for a given embedding dimension remains uncertain.
For fixed accuracy and failure probability, the total cost of \krill preconditioning, including the CG iterations, is therefore $\order((N + k^2)\, k \log k)$ in theory and $\order(N k \log k + k^3)$ in practice.

\subsubsection{Empirical performance}

\begin{figure}[t]    

    \begin{minipage}[b]{0.49\textwidth}
        \centering
        \quad Large regularization: $\mu = 10^{-6} N$
            
        \vspace{0.5em}
        \includegraphics[width = \textwidth]{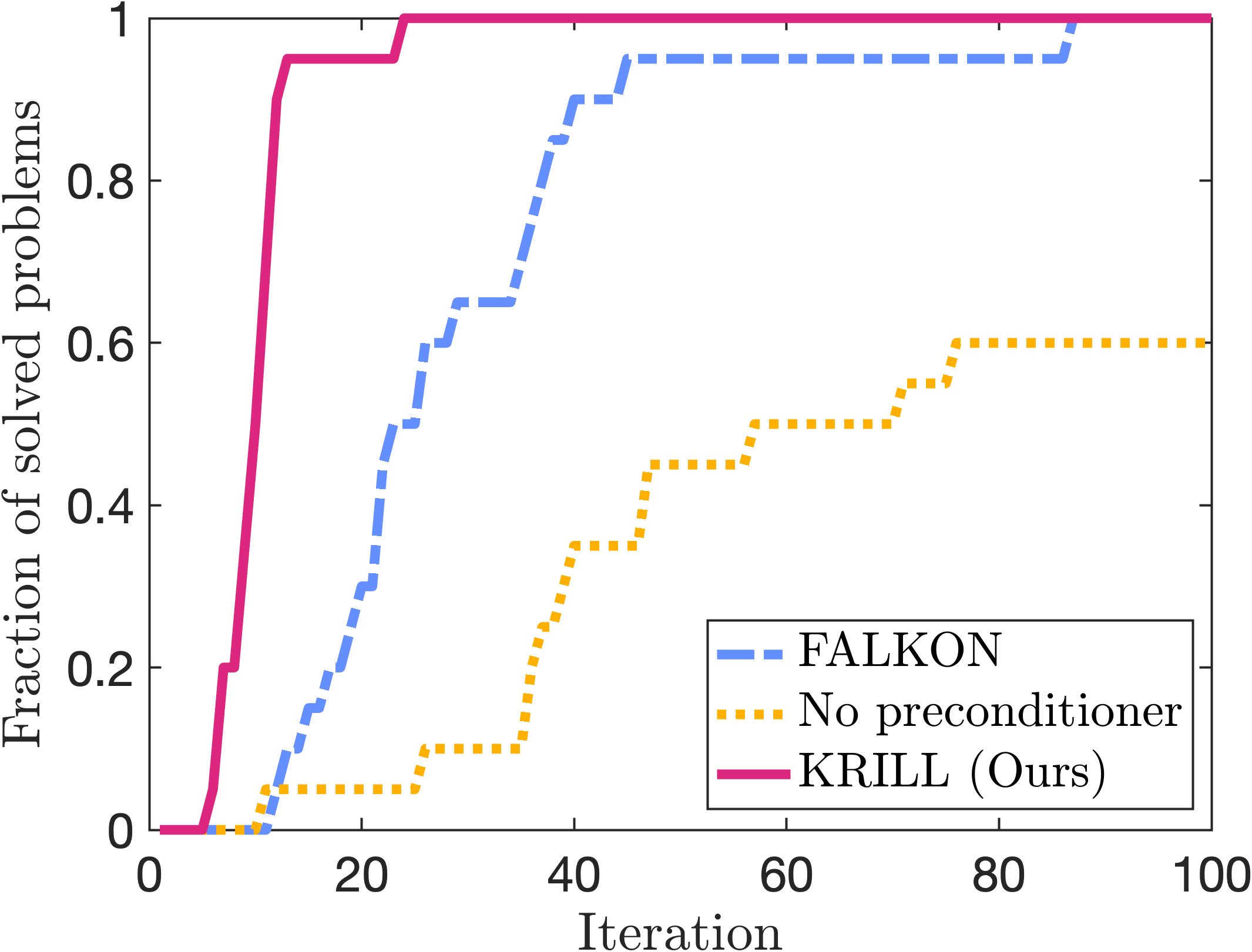}
        \end{minipage}
    \hfill
    \begin{minipage}[b]{0.49\textwidth}
        \centering
        \quad Tiny regularization: $\mu = 10^{-12} N$
        
        \vspace{0.5em}
        
        \includegraphics[width = \textwidth]{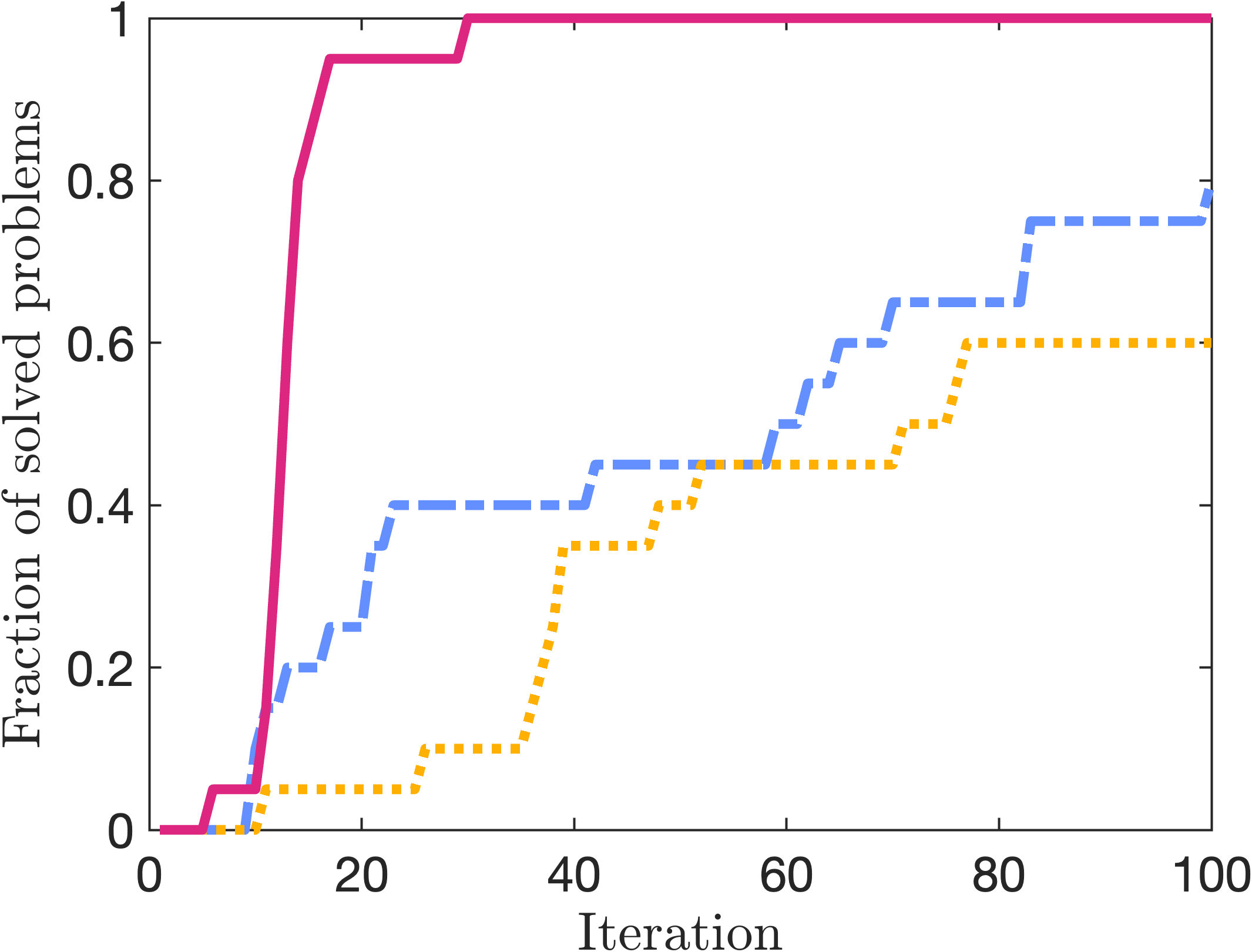}
    \end{minipage}\hfill
    \caption{Fraction of solved problems versus number of CG iterations for the 20 kernel problems in \Cref{tab:datasets}.}
    \label{fig:KRILL-performance}
\end{figure}

\Cref{fig:KRILL-performance} evaluates the performance of \krill across the 20 regression and classification problems described in \Cref{tab:datasets}.
For each problem, we subsample $N = 4 \times 10^4$ data points 
and select $k = 1000$ centers uniformly at random.
We standardize the features and apply KRR using the squared exponential kernel \eqref{eq:squared_exponential} with the bandwidth $\sigma = \sqrt{d}$.
We either set the regularization to a large value $\mu/N = 10^{-6}$ or a tiny value $\mu/N = 10^{-12}$.
We run $100$ CG iterations and declare a problem to be ``solved'' if the relative residual
\begin{equation}
\label{eq:relative_res}
    \frac{\lVert \bigl[\mat{A}(\set{S}, :) \mat{A}(:, \set{S}) + \mu \mat{A}(\set{S}, \set{S})\bigr] \vec{\hat{\beta}} - \mat{A}(\set{S}, :) \vec{y} \rVert}{\lVert \mat{A}(\set{S}, :) \vec{y} \rVert}
\end{equation}
falls below a tolerance of $\varepsilon = 10^{-4}$.
For comparison, \Cref{fig:KRILL-performance} also evaluates the performance of the \FALKON preconditioner that approximates the Gram matrix $\mat{G}$ via Monte Carlo sampling.
More details about \FALKON are in \Cref{sec:restricted_preconditioners}.

Examining \Cref{fig:KRILL-performance}, we find that \krill solves all classification and regression problems in 30 or fewer iterations, both when the regularization is large (left) and when the regularization is tiny (right).
In contrast, \FALKON struggles with the tiny regularization (right panel), and it solves just 8 problems after 30 CG iterations.

\begin{figure}[t]    
    \begin{minipage}[b]{0.49\textwidth}
    \centering
    \quad \texttt{COMET\_MC\_SAMPLE}
    \vspace{0.5em}
    \includegraphics[width = \textwidth]{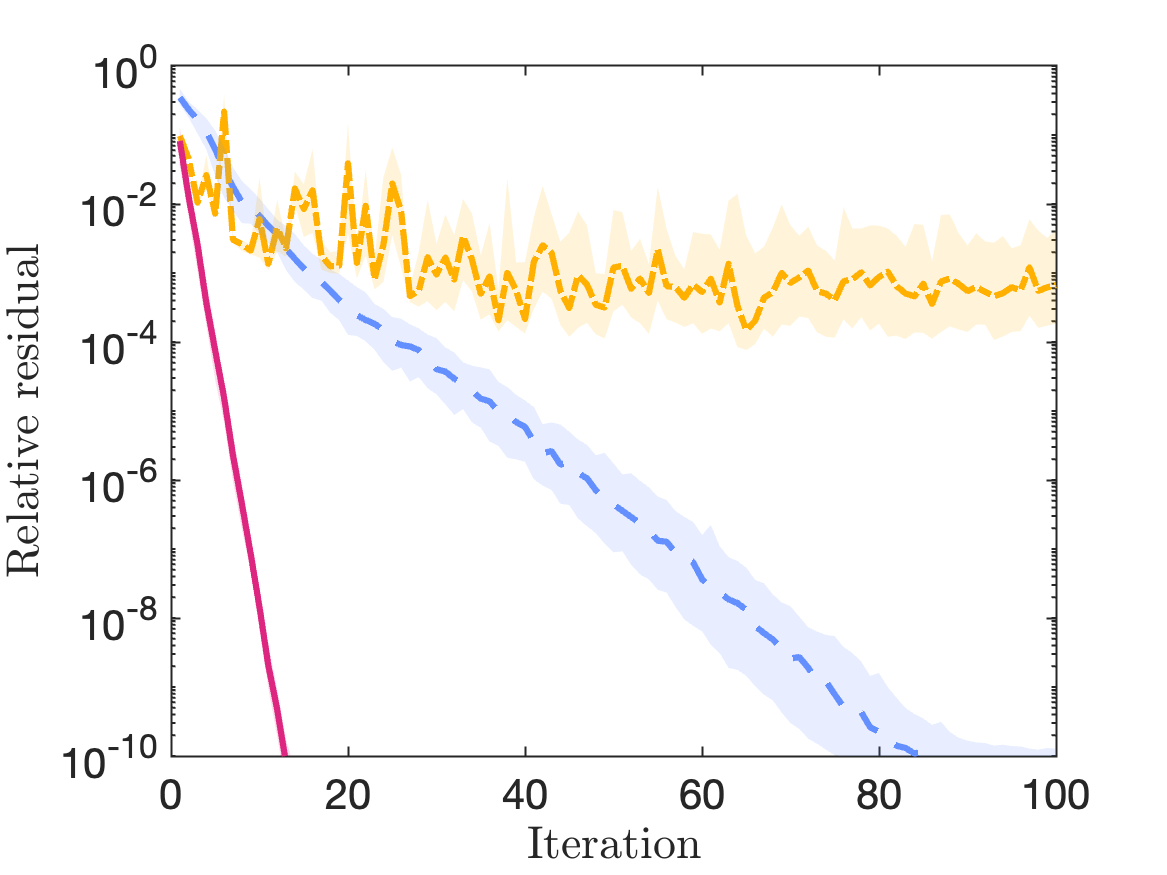}
    \end{minipage}
    \hfill
    \begin{minipage}[b]{0.49\textwidth}
    \centering
        \quad \texttt{creditcard}
    \vspace{0.5em}
    \includegraphics[width = \textwidth]{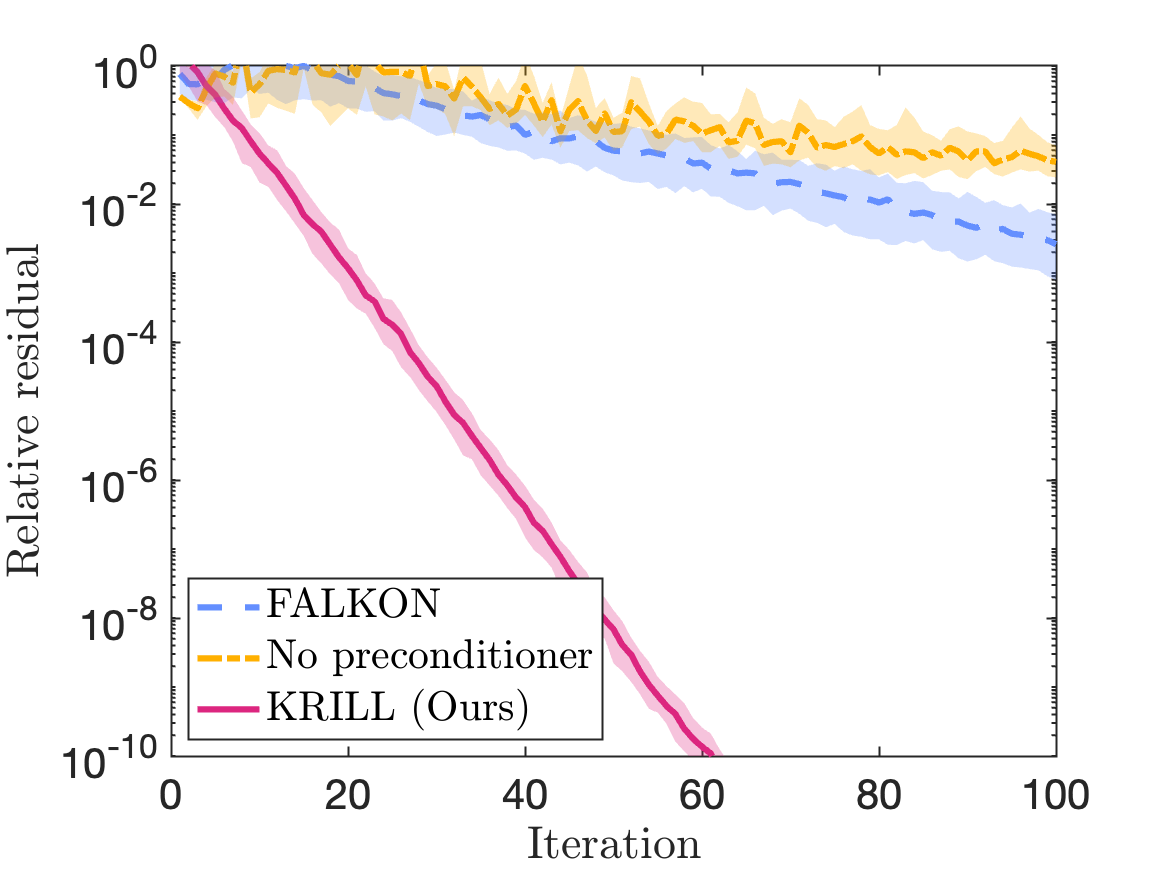}
    \end{minipage}\hfill
    \caption{Relative residual versus number of CG iterations for the problems with the fastest (\texttt{COMET\_MC\_SAMPLE}, left) and slowest (\texttt{creditcard}, right) CG convergence when $\mu/N = 10^{-6}$.
    Thick lines show the median and shaded regions show the $20\%$--$80\%$ error quantiles over 100 random trials.}
    \label{fig:KRILL-HIGGS}
\end{figure}

\krill is a randomized algorithm, so its behavior depends on the random seed.
Yet the results of \krill are consistent up to $\pm 20\%$ random variation in the number of CG iterations required for convergence.
\Cref{fig:KRILL-HIGGS} shows the median and the 20\%--80\% quantiles of the relative residual \eqref{eq:relative_res} for the problems with the fastest (\texttt{COMET\_MC\_SAMPLE}, left) and slowest (\texttt{creditcard}, right) CG convergence.
The quantiles are calculated using 100 independent trials including random subsampling of the data set and random preconditioner formation.
The quantiles are tightly concentrated around the median.

\subsubsection{Theoretical guarantees}

With the proper parameter choices, \krill is guaranteed to solve any nonsingular restricted KRR problem.
\textbf{Unlike \RPCholesky preconditioning, \krill does not require any eigenvalue decay.}
Here, we state our main theoretical guarantee, which is proved in \Cref{sec:a_priori_2}.

\begin{theorem}[\krill performance guarantee]
\label{thm:a_priori_2}
Fix a failure probability $\delta \in (0,1)$ and an error tolerance $\epsilon \in (0,1)$.  
Let $\mat{A}$ be any positive-semidefinite matrix, let $\set{S}$ be an index set of cardinality $|\set{S}| = k$, and let $\mu > 0$ be a positive number.
Define $\mat{M} = \mat{A}(\set{S}, :) \mat{A}(:, \set{S}) + \mu \mat{A}(\set{S}, \set{S})$ and assume $\mat{M}$ is nonsingular.
Generate a sparse random sign embedding $\mat{\Phi} \in \real^{d \times N}$ where the column sparsity $\zeta$ and embedding dimension $d$ are functions of $k$ and $\delta$ given in \cite{Coh16} that satisfy
\begin{equation*}
\zeta = \order(\log(k \slash \delta)) \quad \text{and} \quad d = \order(k \log(k \slash \delta)).
\end{equation*}
With probability at least $1 - \delta$, the \krill preconditioner $\mat{P} = \mat{A}(\set{S}, :) \mat{\Phi}^\ast \mat{\Phi} \mat{A}(:, \set{S}) + \mu \mat{A}(\set{S}, \set{S})$ controls the condition number at a level
\begin{equation} \label{eq:krill-success-cond}
    \kappa\bigl(\mat{P}^{-1/2} \mat{M} \mat{P}^{-1/2}\bigr) \leq 3.
\end{equation}
Conditional on the event~\eqref{eq:krill-success-cond}, when we apply PCG to the linear system $\mat{M} \hat{\vec{\beta}} = \mat{A}(\set{S}, :) \vec{y}$, we obtain an approximation $\hat{\vec{\beta}}^{(t)}$ to the actual solution $\hat{\vec{\beta}}$ that satisfies
\begin{equation*}
    \lVert \hat{\vec{\beta}}^{(t)} - \hat{\vec{\beta}} \rVert_{\mat{M}}
    \leq \epsilon \lVert \hat{\vec{\beta}} \rVert_{\mat{M}},
\end{equation*}
at any iteration
$t \geq \log (2 \slash \epsilon)$.
\end{theorem}

\Cref{thm:a_priori_2} implies that \krill can solve any nonsingular restricted KRR problem in $\order((N + k^2)\, k \log k)$ operations.
In contrast, \FALKON and related preconditioners \cite{RCR17,meanti2020kernel,RCCR19} are only guaranteed to solve restricted KRR problems when the number of centers and the regularization satisfy $k = \Omega(\sqrt{N})$ and $\mu = \Omega(\sqrt{N})$ \cite{RCR17,RCCR19}.
These restrictions help explain why \krill is more robust than \FALKON in our experiments.

\begin{figure}[t]
    \centering
    \quad \qquad Original eigenvalues
    \quad \qquad \krill eigenvalues
    \qquad Eigenvalues $\lambda_i((\mat{Z}^* \mat{Z})^{-1})$

    \qquad \qquad \quad
    $\lambda_i(\mat{M})$
    \qquad \qquad \qquad $\lambda_i(\mat{P}^{-1/2} \mat{M} \mat{P}^{-1/2})$
    \qquad \quad for $\mat{Z} \sim \mathcal{N}(0, 1)^{2k \times k}$

    \vspace{0.5em}
    \includegraphics[width = \textwidth]{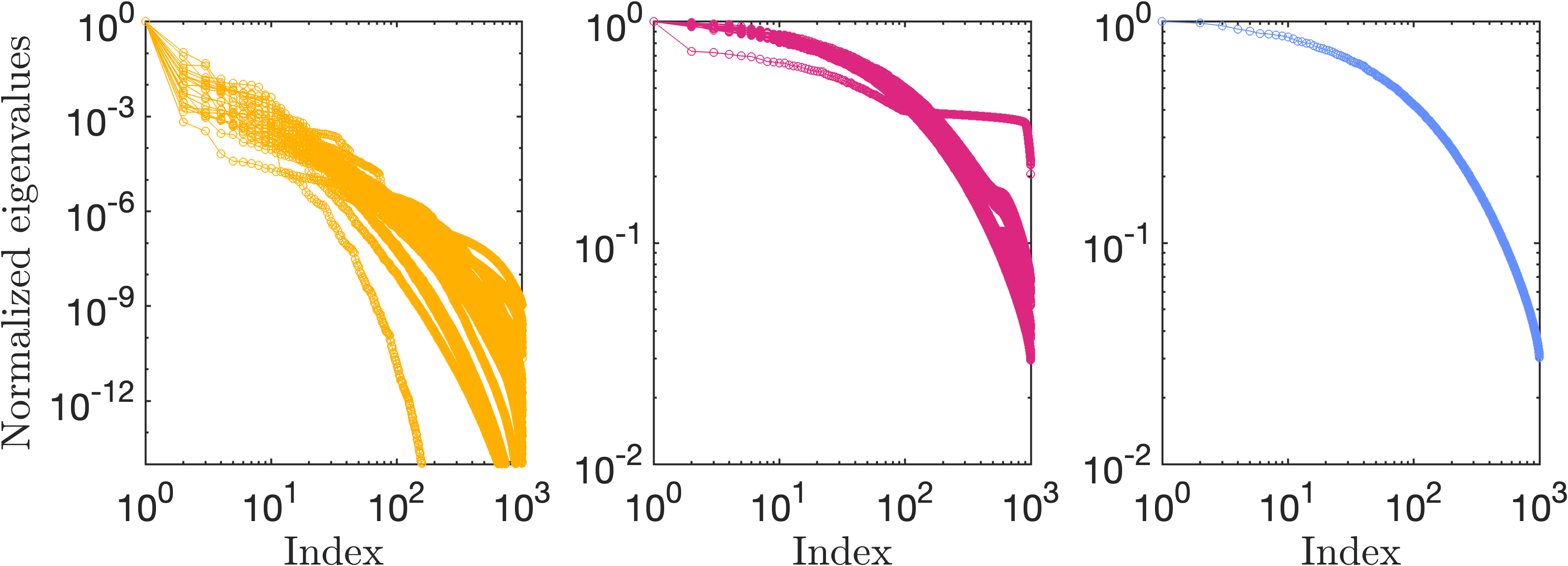}
    \caption{Eigenvalue decay for the 20 KRR problems in \Cref{tab:datasets}.
    Left panel shows eigenvalues of the matrix $\mat{M} = \mat{A}(\set{S}, :) \mat{A}(:, \set{S}) + \mu \mat{A}(\set{S}, \set{S})$ before \krill preconditioning (normalized by the largest eigenvalue).
    Middle panel shows normalized eigenvalues of the matrix $\mat{P}^{-1/2} \mat{M} \mat{P}^{-1/2}$ after \krill preconditioning.
    Right panel shows eigenvalues $\lambda_i((\mat{Z}^* \mat{Z})^{-1})$ for a Gaussian matrix $\mat{Z} \sim \mathcal{N}(0, 1)^{2k \times k}$.
    Note the different vertical axis scales.}
    \label{fig:eigenvalues_approx}
\end{figure}

\Cref{fig:eigenvalues_approx} evaluates \krill's impact on the 20 example problems with a tiny regularization $\mu = 10^{-12} N$.
Before \krill preconditioning, the condition number $\lambda_{\max} / \lambda_{\min}$ ranges from $10^8$ to $10^{15}$ (left panel).
After \krill preconditioning, the preconditioned eigenvalues (middle panel) are closely modeled by the inverse squared singular values of a $2k \times k$ Gaussian matrix (right panel), leading to a condition number $(1 + 1/\sqrt{2})^2 / (1 - 1/\sqrt{2})^2 \approx 34.0$.
The Gaussian comparison was presented in \cite[Sec.~2]{ozaslan2019iterative}, and it represents a broader Gaussian universality phenomenon in random matrix theory \cite{tropp2025comparisontheoremsminimumeigenvalue,chenakkod2025optimalsubspaceembeddingsresolving}.
The comparison suggests that \krill's performance depends weakly on the particular kernel matrix and is governed largely by the embedding aspect ratio $d/k=2$.

\section{Background and comparisons with other preconditioners} \label{sec:history}

In this section, we compare our new conjugate gradient (CG) preconditioners with existing preconditioners for solving the full-data and restricted KRR equations.

\paragraph{Preconditioned conjugate gradient}

Both the full-data and restricted KRR equations have the form
\begin{equation}
\label{eq:psd_system}
    \mat{M}\vec{z} = \vec{b},
\end{equation}
where $\mat{M}$ is positive definite. Given a positive-definite preconditioner $\mat{P}$, we apply standard PCG with a zero initial iterate and terminate when the relative residual falls below $\varepsilon$ \cite[\S9.2]{Saa03}. Equivalently, PCG applies CG to $\mat{P}^{-1/2} \mat{M} \mat{P}^{-1/2}$, so its convergence is controlled by
\begin{equation*}
\kappa_{\mat{P}}=\kappa\bigl(\mat{P}^{-1/2} \mat{M} \mat{P}^{-1/2} \bigr).
\end{equation*}
After $t$ iterations, we can bound the distance between the PCG iterate $\vec{z}^{(t)}$ and the solution $\vec{z}$ to the system \eqref{eq:psd_system} as follows \cite[eq.~6.128]{Saa03}:
\begin{equation}
\label{eq:rate}
    \frac{\lVert \vec{z}^{(t)}-\vec{z}\rVert_{\mat{M}}}{\lVert \vec{z} \rVert_{\mat{M}}}
    \leq
    2\left(\frac{\sqrt{\kappa_{\mat{P}}}-1}
{\sqrt{\kappa_{\mat{P}}}+1}\right)^t.
\end{equation}
Thus, fast convergence occurs when $\mat{P}$ is a good spectral approximation to $\mat{M}$.

In the rest of this section, we review the history of PCG methods for the full-data and restricted KRR problems.
The section does not cover KRR methods based on sketch-and-solve \cite{bach2013sharp,alaoui2015fast} or row access methods like stochastic gradient descent \cite{dai2014scalable,ma2017diving}, 
because our focus is PCG preconditioning.

\subsection{Preconditioners for full-data KRR}
\label{sec:exact_krr_review}

To efficiently solve the full-data KRR problem $(\mat{A} + \mu \Id) \vec{x} = \vec{y}$, 
we can apply a preconditioner of the form
\begin{equation*}
    \mat{P} = \mat{\hat{A}} + \mu \Id,
\end{equation*}
where $\mat{\hat{A}} \in \real^{N \times N}$ is an approximation of the full kernel matrix $\mat{A}$.
Ideally, the approximation $\mat{\hat{A}} \approx \mat{A}$ is highly accurate, and we can quickly apply $\mat{P}^{-1}$ to a vector.

In the KRR literature,
$\mat{\hat{A}}$ is typically constructed as a low-rank approximation \cite{COCF16,avron2017faster,gardner2018gpytorch,wang2019exact,FTU21}.
A low-rank approximation can only be accurate when the kernel matrix is close to being low-rank, i.e., when the eigenvalues of $\mat{A}$ decay quickly.
However, even under favorable eigenvalue decay conditions, it has been a challenge to identify accurate and computationally tractable low-rank approximations.

One low-rank approximation used in KRR problems is the \emph{column Nystr\"om approximation} \cite[\S19.2]{MT20a}, which is
obtained from a partial Cholesky decomposition with pivoting.
The column Nystr\"om approximation takes the form
\begin{equation*}
    \mat{\hat{A}} = \mat{A}(:, \textsf{S}) \mat{A}(\textsf{S}, \textsf{S})^\dagger \mat{A}(\textsf{S}, :),
\end{equation*}
where $\textsf{S} = \{s_1, s_2, \ldots, s_r\}$ is the set of pivots chosen in the Cholesky procedure.
The column Nystr\"om approximation is convenient because it is formed from the columns of $\mat{A}$ indexed by $\textsf{S}$ and does not require viewing the rest of the matrix.
Yet the approximation accuracy depends on the index set $\set{S}$.
In the context of full-data KRR, the most popular strategies for picking the set $\set{S}$ of columns are uniform sampling and greedy selection.
Both of these approaches exhibit failure modes.

\paragraph{Failure of uniform sampling}
One simple strategy is to select the column pivots $s_1, \ldots, s_r$ uniformly at random \cite{COCF16,FTU21}.
This \emph{uniform sampling} method can be effective for some problems, but it fails to explore less populated regions of data space.
To illustrate this shortcoming, let $\vec{0}_m$ and $\vec{1}_m$ denote the $m$-vector with entries equal to $0$ and $1$, respectively, and consider the kernel matrix
\begin{equation*}
        \mat{A} = \begin{pmatrix}
        \vec{1}_{N-N^{1/3}}^{\vphantom{*}} \vec{1}_{N-N^{1/3}}^\ast & \vec{0}_{N-N^{1/3}}^{\vphantom{*}} \vec{0}_{N^{1/3}}^* \\
        \vec{0}_{N^{1/3}}^{\vphantom{*}} \vec{0}_{N-N^{1/3}}^* & \vec{1}_{N^{1/3}}^{\vphantom{*}} \vec{1}_{N^{1/3}}^\ast
    \end{pmatrix}.
\end{equation*}
Constructing a low-rank approximation for $\mat{A}$ should be easy,
since the rank is $2$.
However, uniform sampling selects many columns from the left block and neglects columns from the right block that are also needed to ensure the approximation quality.
Uniform sampling needs an approximation rank much higher than $2$ and at least as large as $r = \Omega(N^{2/3})$ to build an effective preconditioner.

\paragraph{Failure of greedy selection} A second strategy for choosing columns is the \emph{greedy selection} method \cite{gardner2018gpytorch,wang2019exact},
in which we adaptively select each column pivot $s_i$ by finding the largest diagonal element of the residual matrix $\mat{A} - \mat{\hat{A}}_{(i-1)}$ after $i-1$ steps of the Cholesky procedure.
The greedy method has the opposite failure mode from the uniform method: it focuses on outlier data points and fails to explore highly populated regions of data space.
To see this limitation, let $\Id_{N^{2/3}}$ denote the $N^{2/3} \times N^{2/3}$ identity matrix, and consider the kernel matrix
\begin{equation*}
    \mat{A} = \vec{1}_N^{\vphantom{\ast}} \vec{1}_N^\ast + 
    \begin{pmatrix}
        \frac{1}{2N} \, \vec{1}_{N-N^{2/3}}^{\vphantom{*}} \vec{1}_{N-N^{2/3}}^\ast & \vec{0}_{N-N^{2/3}}^{\vphantom{*}} \vec{0}_{N^{2/3}}^* \\
        \vec{0}_{N^{2/3}}^{\vphantom{*}}
        \vec{0}_{N - N^{2/3}}^*
        & \frac{1}{N} \, \Id_{N^{2/3}}
    \end{pmatrix}.
\end{equation*}
Constructing a preconditioner for $\mat{A}$ should be easy, since the $\mu$-tail rank (\Cref{def:d_tail}) is $\textup{rank}_{\mu}(\mat{A}) \leq 2$ for any $\mu \geq N^{-1/3}$.
Yet the greedy strategy selects columns from the right block first and misses the left columns.
To provide an effective preconditioner,
greedy sampling needs an approximation rank $r > N^{2/3}$, which is large enough that all the columns in the right block have been selected.

\paragraph{Advantages of \RPCholesky} In this paper, we propose a new KRR preconditioner that uses \RPCholesky (see \Cref{alg:rpcholesky} or \cite{chen2022randomly}) to select the columns for the Nystr\"om approximation.
In \RPCholesky, we adaptively sample the column pivot $s_i$ with probability proportional to the diagonal entries of the residual matrix $\mat{A} - \mat{\hat{A}}_{(i-1)}$ after $i-1$ steps of the Cholesky procedure.
Because of the adaptive sampling distribution, \RPCholesky balances exploration of the small diagonal entries and exploitation of the large diagonal entries, avoiding the failure modes of both the uniform and greedy strategies.
\RPCholesky is guaranteed to produce an effective preconditioner if the approximation rank $r$ is set to a modest multiple of the $\mu$-tail rank (\Cref{thm:a_priori}).

\paragraph{Ridge leverage score sampling}
\emph{Ridge leverage score sampling} \cite{alaoui2015fast,musco2017recursive,RCCR19} is a more elaborate procedure for selecting the pivots in the Nystr\"om approximation.
This method first approximates the ridge leverage scores, which are defined as the diagonal entries $\vec{\ell}^{\lambda} = \diag(\mat{A}(\mat{A} + \lambda \Id)^{-1})$ for a parameter $\lambda > 0$.
Then, the method randomly generates a Nystr\"om approximation whose pivots are selected with probabilities proportional to the ridge leverage scores.
While RLS sampling can lead to an approximation that is nearly as accurate as \RPCholesky (see \Cref{fig:exact-performance}), it can also produce lower-quality approximations given challenging inputs \cite{chen2022randomly}.
Additionally, RLS sampling is slower than the blocked implementation of \RPCholesky given in \Cref{alg:rpcholesky}, because of the requirement to approximate the ridge leverage scores.
For example, when we apply RLS sampling with the RecursiveRLS algorithm \cite{musco2017recursive} to compute a rank-$1000$ approximation of a kernel matrix of size $N = 4 \times 10^4$, the runtime is $2.3\times$ slower than \RPCholesky's runtime.

\paragraph{More expensive Gaussian preconditioner} A different type of preconditioner for full-data KRR is based on the \emph{Gaussian Nystr\"om approximation} \cite{FTU21}
\begin{equation*}
    \mat{\hat{A}} = \mat{A}\mat{Z} (\mat{Z}^*\mat{A}\mat{Z})^\dagger \mat{Z}^*\mat{A},
\end{equation*}
where $\mat{Z} \in \real^{N\times r}$ is a matrix with independent standard normal entries.
Forming the
Gaussian Nystr\"om approximation 
requires $r$ matrix--vector multiplications with the full kernel matrix. 
%which is much more expensive than forming a column Nystr\"om approximation \eqref{eq:column_nys} in factored form.
To maintain a cost of $\order(N^2)$ operations, we can only run Gaussian Nystr\"om with a constant approximation rank $r = \order(1)$, whereas we can run \RPCholesky preconditioning with a much larger approximation rank $r = \order(\sqrt{N})$.
Because of the higher approximation rank,
\RPCholesky preconditioning leads to a more accurate approximation with stronger guarantees.

\paragraph{Other approaches} We briefly mention two other types of approximations that can be used to build preconditioners for the full-data KRR problem.
First, certain types of kernel matrices can be approximated using \emph{random features} \cite{COCF16,avron2017faster}, but the quality of the approximation tends to be poor.
Numerical tests indicate that the random features approach does not yield a competitive preconditioner~\cite{COCF16,YLM+12}.
See \Cref{fig:exact-performance,fig:slower-methods} for further numerical comparisons.

Second, kernel matrices can be approximated using carefully constructed hierarchical approximations \cite{COCF16,AFG+16,CAS17,YLRB17}
or sparse inverse approximations \cite{schafer2021sparse,zhao2024adaptive,huan2025sparse}.
These approaches extend beyond low-rank structure and provide richer matrix approximation classes; however, the applicability to high-dimensional data sets remains unclear \cite{AFG+16}.
The strongest existing guarantees and empirical demonstrations for these approaches generally concern low-dimensional input data.

\subsection{Preconditioners for restricted KRR} \label{sec:restricted_preconditioners}

To efficiently solve the restricted KRR problem $\bigl[\mat{A}(\set{S}, :) \mat{A}(:, \set{S}) + \mu \mat{A}(\set{S}, \set{S})\bigr]\, \vec{\hat{\beta}} = \mat{A}(\set{S}, :) \vec{y}$,
we can apply a preconditioner of the form
\begin{equation*}
    \mat{P} = \mat{\hat{G}} + \mu \mat{A}(\set{S}, \set{S}),
\end{equation*}
where $\mat{\hat{G}} \in \real^{k \times k}$ approximates the Gram matrix $\mat{G} = \mat{A}(\set{S}, :) \mat{A}(:, \set{S})$.

\paragraph{\FALKON-type preconditioners} \FALKON-type preconditioners \cite{RCR17,RCCR19,meanti2020kernel} are based on a Monte Carlo approximation of $\mat{G}$.
The Monte Carlo approach assumes that each data point $\vec{x}^{(i)}$ is selected as a center with nonzero probability $p_i = \mathbb{P}\{i \in \set{S}\}$.
The Gram matrix is then approximated as
\begin{equation}
\label{eq:more_general}
    \mat{\hat{G}} = \mat{A}(\set{S}, \set{S}) \mat{D}(\set{S}, \set{S})^{-1} \mat{A}(\set{S}, \set{S}),
\end{equation}
where $\mat{D} = \diag(p_1, p_2, \ldots, p_N)$ is the diagonal matrix of the selection probabilities.
The matrix entries $g_{ij}$ and $\hat{g}_{ij}$ can be explicitly written as
\begin{equation*}
    g_{ij}
    = \sum_{\ell = 1}^N a_{s_i,\ell} a_{\ell, s_j}
    \quad \text{and} \quad
    \hat{g}_{ij} = \sum_{\ell=1}^N \frac{\mathbf{1}\{\ell \in \set{S}\}}{p_\ell} a_{s_i,\ell} a_{\ell, s_j}.
\end{equation*}
These entries are the same in expectation:
\begin{equation*}
    \sum_{\ell = 1}^N a_{i,\ell} a_{\ell, j} 
    = \expect\Biggl[ \sum_{\ell=1}^N \frac{\mathbf{1}\{\ell \in \set{S}\}}{p_\ell} a_{i,\ell} a_{\ell, j} \Biggr] \quad \text{for fixed indices $1 \leq i,j \leq N$,}
\end{equation*}
which suggests that the quantities $g_{ij}$ and $\hat{g}_{ij}$ will be close under appropriate conditions. 
\FALKON \cite{RCR17,meanti2020kernel} is a special case of \eqref{eq:more_general} in which the centers are sampled uniformly at random, and \FALKON-\BLESS \cite{RCCR19} is a special case of \eqref{eq:more_general} in which the centers are chosen by ridge leverage score sampling.

The main limitation of \FALKON-type preconditioners is the Monte Carlo approximation error.
To account for the Monte Carlo error, the available analyses \cite{RCR17,RCCR19} assume a large number of centers $k = \Omega(\sqrt{N})$ and a large regularization $\mu = \Omega(\sqrt{N})$.
When these assumptions are satisfied, \FALKON-type preconditioners solve the restricted KRR equations in $\order((N + k^2) k)$ operations.
However, the assumptions on $k$ and $\mu$ may not be satisfied in practice.
For example, $\mu$ is often small when chosen by cross-validation, and \cite{RCR17} contains examples of $\mu$ values up to five orders of magnitude smaller than $\mu = \sqrt{N}$.

\paragraph{\krill} In this paper, we recommend the approximation
\begin{equation*}
\mat{\hat{G}} = \mat{A}(\set{S}, :) \mat{\Phi}^\ast \mat{\Phi} \mat{A}(:, \set{S}),
\end{equation*}
where $\mat{\Phi} \in \real^{d \times N}$ is the sparse sign embedding defined in \cref{sec:krill_construction}.
Whereas the \FALKON preconditioner \eqref{eq:more_general} forms a preconditioner from the square submatrix $\mat{A}(\set{S}, \set{S})$, \krill uses the full matrix $\mat{A}(:,\set{S})$, opening up the possibility for higher accuracy.
We can guarantee the effectiveness of the \krill preconditioner 
$\mat{P} = \mat{\hat{G}} + \mu \mat{A}(\set{S}, \set{S})$,
regardless of the number of centers $k$ and the regularization $\mu$. 

\krill is inspired by the \emph{sketch-and-precondition} approach for solving overdetermined least-squares problems \cite[\S10.5]{MT20a}.
In sketch-and-precondition, we use a random embedding to approximate the matrix appearing in the normal equations, 
and the approximation serves as a preconditioner for solving the least-squares problem.
Sketch-and-precondition was first proposed in \cite{RT08} and later refined in \cite{AMT10,meng2014lsrn,clarkson2017low,lacotte2021fast,meier2022randomized}.
It has been applied with several different embeddings \cite[\S\S8--9]{MT20a},
but the empirical comparisons of \cite[Fig.~1]{DM21} suggest the sparse sign embedding is most efficient.

\krill differs from classical presentations of sketch-and-precondition because the embedding is applied to the Gram matrix term $\mat{A}(\set{S},:)\mat{A}(:,\set{S})$ but not the regularization term $\mu \mat{A}(\set{S},\set{S})$.
This idea of ``partial sketching'' appears in \cite{pilanci2017newton,murray2022}, and it improves the accuracy of the approximation.
Notwithstanding this difference,
\krill is motivated by the same ideas as classical sketch-and-precondition and its analysis follows the same patterns.

\section{Case studies} \label{sec:case_studies}

In this section, 
we evaluate \RPCholesky on a molecular-regression problem and \krill on a particle-classification problem.

\subsection{HOMO energy prediction}
\label{sec:homo_energy}

The QM9 data set describes the properties of $1.3\times 10^5$ organic molecules \cite{DBB21,WM21}. Following \cite{STR+19}, we predict the highest-occupied-molecular-orbital (HOMO) energy from standardized Coulomb-matrix features using the $\ell_1$ Laplace kernel
\begin{equation*}
    K(\vec{x}, \vec{x}') = \exp \bigl( -\tfrac{1}{\sigma} \lVert \vec{x} - \vec{x}' \rVert_{
\ell_1}\bigr).
\end{equation*}
The bandwidth $\sigma$ and ridge parameter $\mu$ are chosen using cross-validation.
Because the complete kernel matrix exceeds laptop-scale memory, the experiments in \cite{STR+19} used random subsets containing at most $32{,}000$ molecules.
At that sample size, however, the predictive accuracy had not yet saturated, motivating our larger-scale experiments.

\begin{figure}[t]
    \centering
    \includegraphics[width = .5\textwidth]{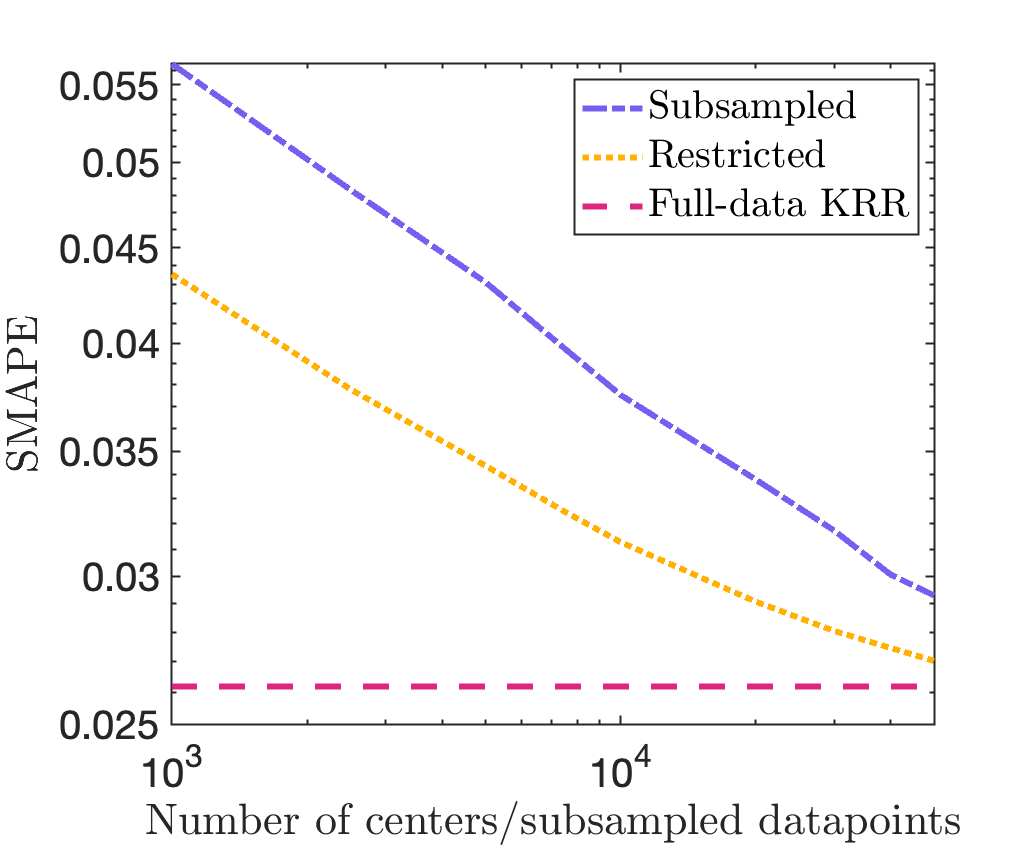}
    \caption{SMAPE versus number of centers / number of subsampled data points for the HOMO energy prediction problem $(N = 10^5)$. }
    \label{fig:accuracy}
\end{figure}

\begin{figure}[t]
    \begin{minipage}[b]{0.49\textwidth}
    \centering
    \includegraphics[width = \textwidth]{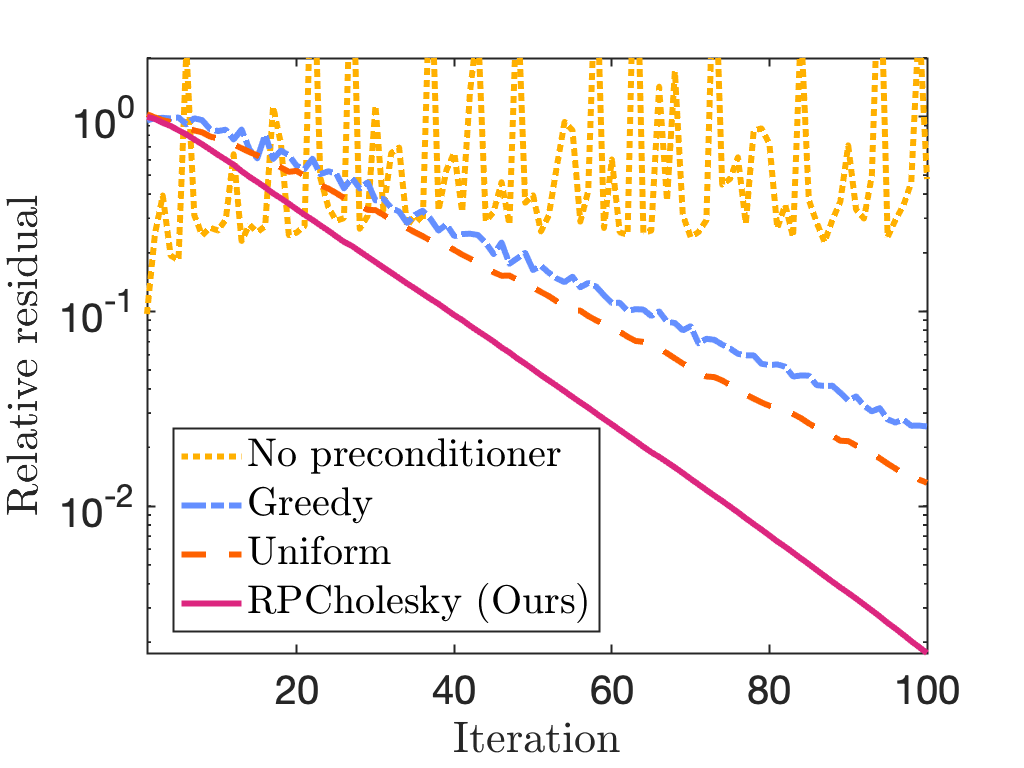}
    \end{minipage}
    \hfill
    \begin{minipage}[b]{0.49\textwidth}
    \centering
    \includegraphics[width = \textwidth]{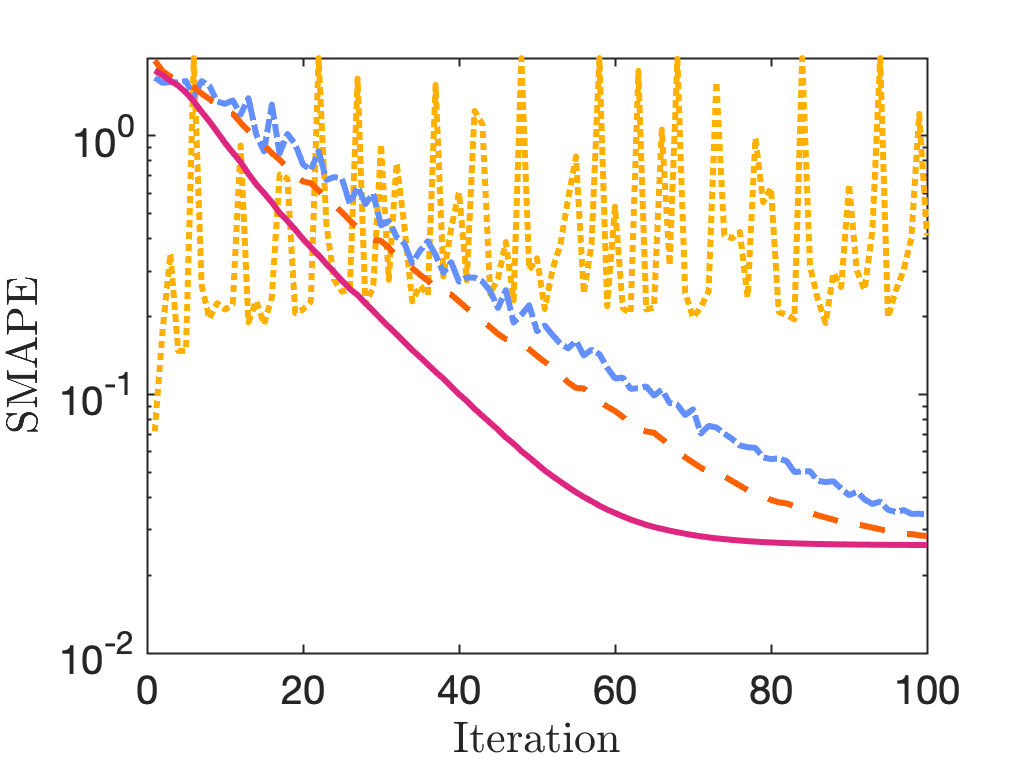}
    \end{minipage}
    \caption{(\emph{left}) Relative residual $\norm{(\mat{A} + \mu \Id)\vec{\beta} - \vec{y}}/\norm{\vec{y}}$ and (\emph{right}) SMAPE  for different column Nystr\"om preconditioners
    for the HOMO energy prediction task.}
    \label{fig:chemistry}
\end{figure}

\begin{figure}[t]
    \centering
    \includegraphics[width = .49\textwidth]{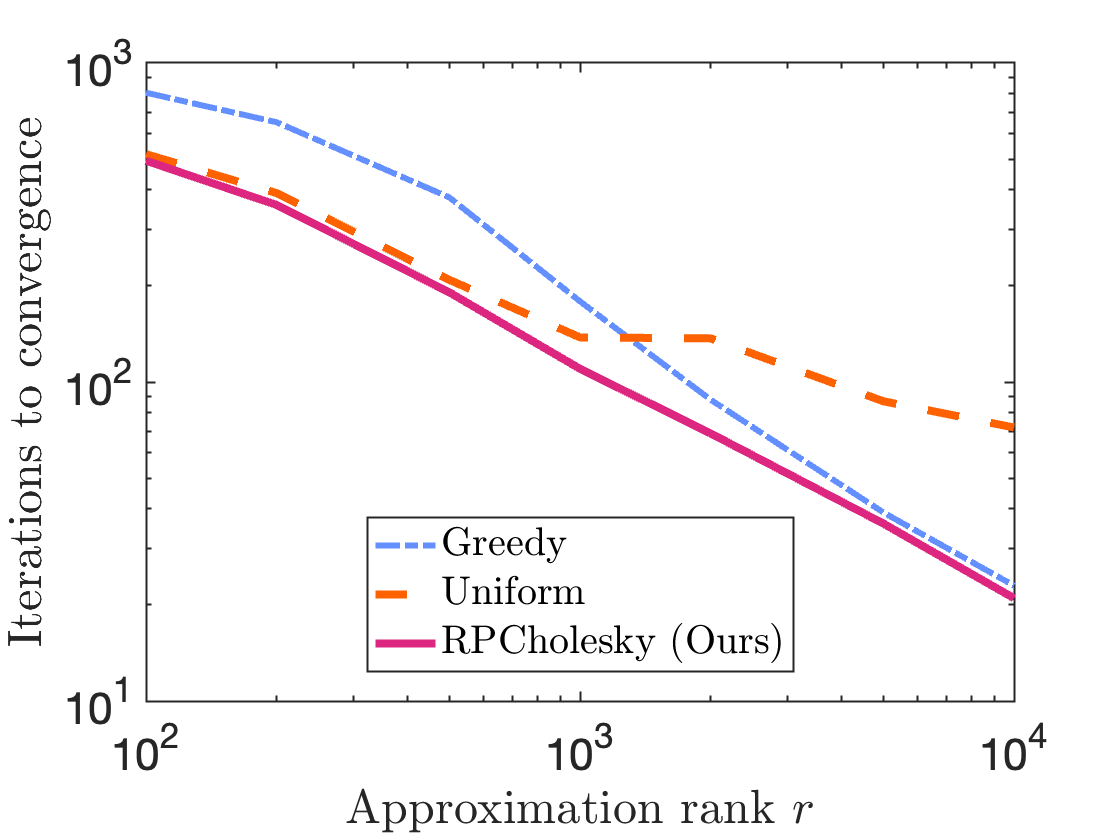}
    \caption{Number of CG iterations to achieve a relative residual of $\varepsilon < 10^{-3}$ with \RPCholesky preconditioning for the HOMO energy prediction task.
    \label{fig:convergence_vs_rank}}
\end{figure}

We can apply KRR to \emph{the complete QM9 data set} either using (i) restricted KRR or (ii) full-data KRR with \RPCholesky preconditioning. 
In all our experiments, we split the data into $N = 10^5$ training points and $M = 3 \times 10^4$ test points, and we measure the predictive accuracy using the symmetric mean absolute percentage error (SMAPE):
\begin{equation*}
    \operatorname{SMAPE}(\vec{\widehat y}, \vec{y}) = \frac{1}{M}\sum_{i=1}^M\frac{|\widehat{y}_i - y_i| }{(|\widehat{y}_i|+|y_i|)/2}.
\end{equation*}
The results are presented in \Cref{fig:accuracy,fig:chemistry,fig:convergence_vs_rank}.

Restricted KRR is a relatively cheap approach for predicting HOMO energies with the QM9 data set, since we can store the $k \times k$ preconditioner and $N \times k$ kernel submatrix in working memory with up to $k = 50{,}000$ centers.
\Cref{fig:accuracy} shows that restricted KRR is more accurate than randomly subsampling the data and applying full-data KRR to the random sample.
Still, the predictive accuracy increases with the number of centers $k$ and does not saturate even when $k = N / 2$.

Full-data KRR is the most accurate approach for predicting HOMO energies, but it relies on repeated matrix--vector multiplications with the full kernel matrix.
Since the kernel matrix is too large to store in 64GB of working memory,
we can only store one block of the matrix at a time and perform the multiplications in a blockwise fashion.
Each multiplication requires tens of minutes of computing, due to the cost of evaluating the kernel matrix each time the entries are needed.
The training of the full-data KRR model is very slow unless we find a preconditioner that controls the total number of CG iterations.

\Cref{fig:chemistry} compares several preconditioning strategies for full-data KRR.
The worst strategy is unpreconditioned CG, which leads to no discernible convergence over the first $100$ iterations.
Uniform and greedy Nystr\"om preconditioning (with an approximation rank $r = 10^3$) improve the convergence speed,
but they still require $100$ or more iterations for the predictive accuracy to saturate.
\RPCholesky (also with $r = 10^3$) is the fastest method, leading to converged accuracy after $60$ iterations.

\Cref{fig:convergence_vs_rank} shows that we can reduce the total number of CG iterations even further, by a factor of five, if we increase the \RPCholesky approximation rank from $r = 10^3$ to $r = 10^4$.
However, we also need to consider the cost of preparing the preconditioner, which requires $\order(N r^2)$ operations.
To achieve an appropriate balance, we recommend tuning the value of $r$ based on computational constraints, and $r$ should typically be small enough that the preconditioner can be stored in working memory.
For example, an intermediate value of $10^3 < r < 10^4$ would be most appropriate for the HOMO energy prediction problem given our computing architecture.

\Cref{fig:convergence_vs_rank} also reveals that \RPCholesky performs similarly to uniform pivot selection for a small approximation rank ($r < 500$) and performs similarly to greedy pivot selection for a large approximation rank ($r > 5000$).
It remains an open problem to identify the problems and approximation ranks for which uniform or greedy pivot selection are as successful as \RPCholesky.

\subsection{Exotic particle detection}
\label{sec:supersymmetry}

Following \cite{RCR17,RCCR19,LLR+22}, we evaluate restricted KRR on the SUSY data set, which contains $5\times 10^6$ simulated particle-collision examples labeled as standard-model or supersymmetric decay processes \cite{baldi2014searching}.

\begin{figure}[t]
    \begin{minipage}[b]{0.49\textwidth}
    \centering
    \includegraphics[width = \textwidth]{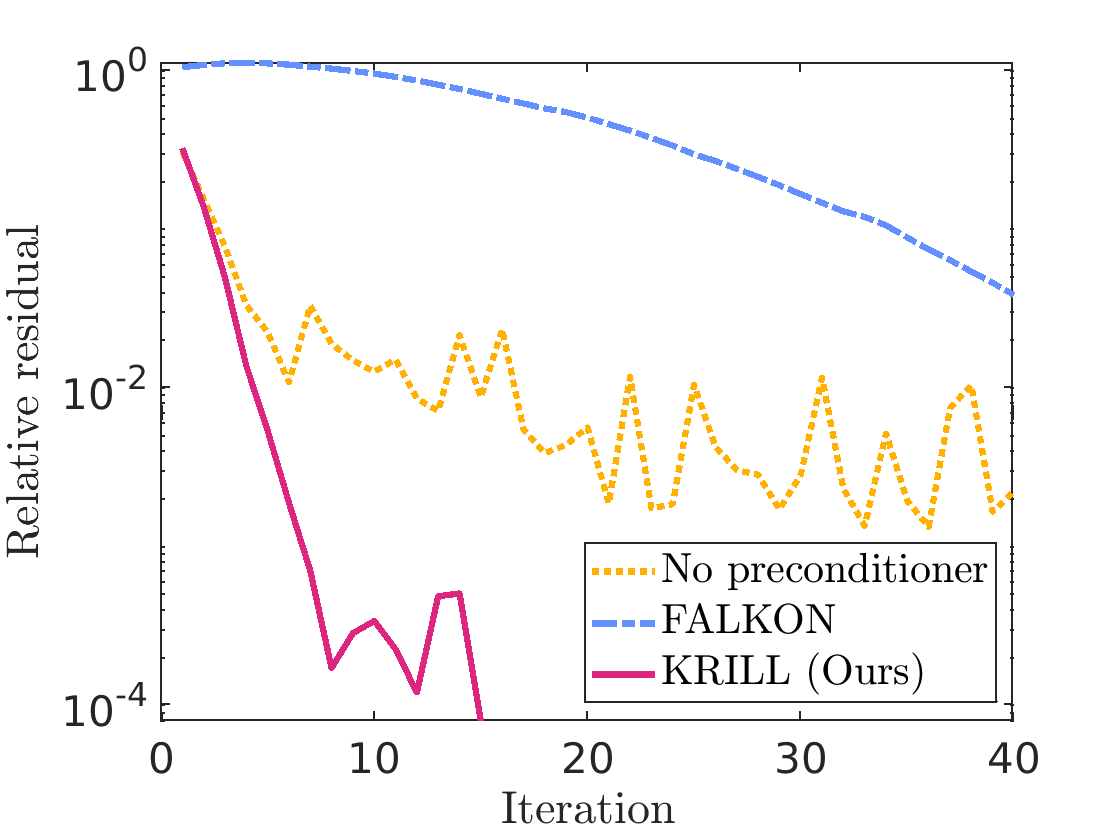}
    \end{minipage}
    \hfill
    \begin{minipage}[b]{0.49\textwidth}
    \centering
    \includegraphics[width = \textwidth]{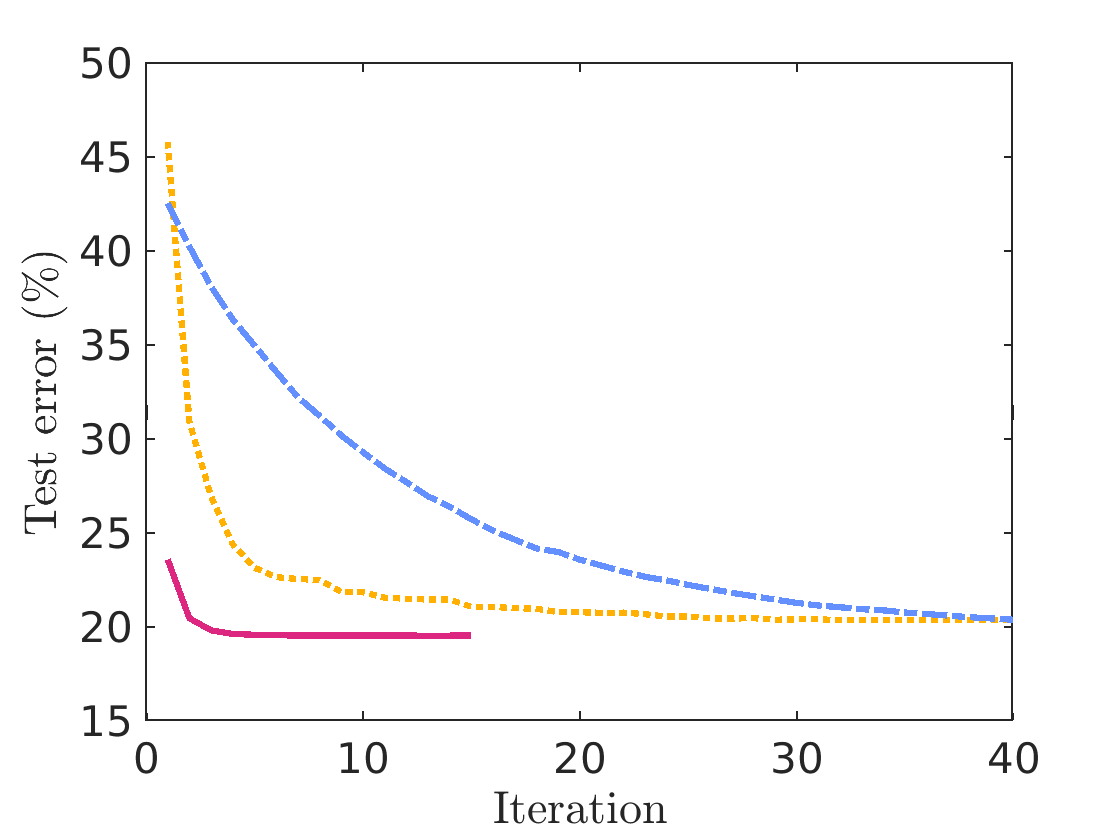}
    \end{minipage}
    \caption{ (\emph{left}) Relative residual $\norm{\smash{ ( \mat{A}(\set{S},:)\mat{A}(:,\set{S}) + \mu \mat{A}(\set{S},\set{S}) )\,\vec{\hat{\beta}}-\mat{A}(\set{S},:)\vec{y}}}/\norm{\mat{A}(\set{S},:)\vec{y}}$ and (\emph{right}) test error for \FALKON vs.\ \krill for the SUSY classification problem.}
    \label{fig:SUSY}
\end{figure}

\begin{figure}[t]
    \centering
    \includegraphics[width=.5\linewidth]{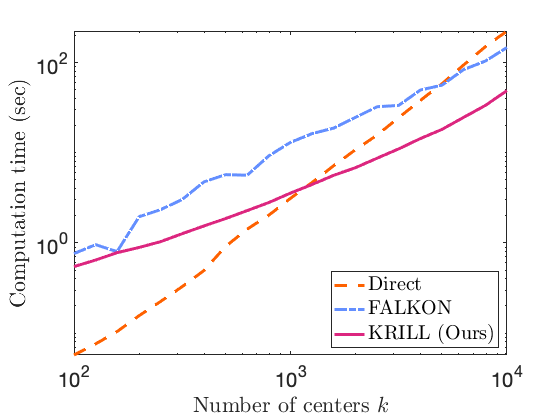}
    \caption{Runtime versus number of centers for the SUSY classification problem ($N = 5 \times 10^5$).}
    \label{fig:lines_crossing}
\end{figure}

For our KRR application, we split the data into $N = 4.5 \times 10^6$ training data points and $M = 5 \times 10^5$ testing data points,
and we select $k = 10^4$ data centers uniformly at random.
Then we standardize the data features
and apply a squared exponential kernel \eqref{eq:squared_exponential} with bandwidth $\sigma = 4$.
The parameters $k$ and $\sigma$ are the same ones used in previous work \cite{RCR17,RCCR19},
but 
we use a smaller regularization $\mu/N = 2\times 10^{-10}$ (compared to $\mu / N = 10^{-6}$ in \cite{RCR17,RCCR19,LLR+22}) which improves the test error from $19.6\%$ to $19.5\%$.
We apply \krill, \FALKON, and unpreconditioned CG to solve the restricted KRR equations with these parameter choices.
We terminate after forty iterations or when the relative residual falls below $\varepsilon = 10^{-4}$, leading to the results in 
\Cref{fig:SUSY,fig:lines_crossing}.

\Cref{fig:SUSY} demonstrates the superior performance of \krill, which reaches a test error of $19.5\%$ after just four iterations.
Meanwhile, unpreconditioned CG fails to achieve such a low test error after forty iterations, and \FALKON reduces the test error even more slowly.
The slow convergence of \FALKON is consistent with the very small normalized regularization $\mu/N=2\times10^{-10}$.The \FALKON preconditioner uses a Monte Carlo approximation based on uniformly sampled Nystr{\"o}m centers.
Under worst-case statistical assumptions, the \FALKON analysis chooses $\mu=\sqrt N$ and requires $k=\Omega(\sqrt N)$ centers \cite[Theorem~3]{RCR17}.
Our parameter choice lies far outside this regime, so the analysis provides no guarantee that the preconditioner improves the spectrum.
For comparison, an earlier \FALKON experiment on SUSY with the larger value $\mu/N=10^{-6}$ attained comparable accuracy after approximately 20 iterations \cite{RCCR19}.
A major advantage of \krill is the theoretically guaranteed performance for all $\mu$ values, which allows users to choose the regularization that gives the lowest test error.

Finally, \Cref{fig:lines_crossing} evaluates the total computation time needed to form the preconditioner and solve the restricted KRR equations up to a relative residual of $\varepsilon = 10^{-4}$.
For this figure, we randomly subsample $N = 5 \times 10^5$ data points so the kernel submatrix fits in 64 GB of working memory, and we compare \krill, \FALKON, and the direct method of forming and factorizing  $\mat{A}(\set{S}, :) \mat{A}(:, \set{S}) + \mu \mat{A}(\set{S}, \set{S})$.
We find that \krill is the dominant algorithm once the number of centers reaches $k = 1250$.
\krill outperforms the direct method because of the superior $\order(Nk\log k + k^3)$ computational cost, and \krill outperforms \FALKON because it converges in roughly $4\times$ fewer CG iterations.
As an advantage, the direct method solves the restricted KRR equations up to machine precision, whereas \krill and \FALKON only solve the equations up to a relative residual of $\varepsilon = 10^{-4}$.
Nonetheless, the test error plateaus at a threshold of $\varepsilon = 10^{-4}$ in all our experiments.

\section{Theoretical results} \label{sec:theory}

In this section, we prove our main theoretical results, \Cref{thm:a_priori,thm:a_priori_2}.

\subsection{Proof of \texorpdfstring{\Cref{thm:a_priori}}{Theorem 1.2}} \label{sec:a_priori}

Let $\mat{\hat{A}}$ be a Nystr{\"o}m approximation of the psd matrix $\mat{A}$, and set $\mat{P} = \mat{\hat{A}} + \mu \Id$.
As the first step of our proof, we will show that
\begin{equation}
\label{eq:condition}
    \kappa \bigl(\mat{P}^{-1 \slash 2} \bigl(\mat{A} + \mu \Id\bigr)\mat{P}^{-1 \slash 2}\bigr) \leq 
    1 + \tr\bigl(\mat{A} - \mat{\hat{A}}\bigr) / \mu.
\end{equation}
To that end, observe that $\mat{P} \preceq \mat{A} + \mu \Id$,
because the Nystr\"om approximation $\mat{\hat{A}}$ is bounded from above by $\mat{A}$.
By conjugation with $\mat{P}^{-1 \slash 2}$, we obtain 
\begin{equation*}
\Id \preceq \mat{P}^{-1 \slash 2} (\mat{A} + \mu \Id) \mat{P}^{-1 \slash 2}
\end{equation*}
and thus
\begin{equation}
\label{eq:conclude}
    1 \leq \lambda_{\rm min} \bigl(\mat{P}^{-1 \slash 2} (\mat{A} + \mu \Id)\mat{P}^{-1 \slash 2}\bigr).
\end{equation}
Next, a short calculation shows that
\begin{equation*}
    \mat{P}^{-1 \slash 2} (\mat{A} + \mu \Id)\mat{P}^{-1 \slash 2}
    = \Id + \mat{P}^{-1 \slash 2} \big(\mat{A} - \mat{\hat{A}}\big)\mat{P}^{-1 \slash 2}.
\end{equation*}
Since the spectral norm is submultiplicative and it equals the largest eigenvalue of a psd matrix, we deduce that
\begin{align*}
    \lambda_{\rm max}\bigl(\mat{P}^{-1 \slash 2} (\mat{A} + \mu \Id)\mat{P}^{-1 \slash 2}\bigr)
    &= 1 + \lambda_{\rm max}\bigl(\mat{P}^{-1 \slash 2} \big(\mat{A} - \mat{\hat{A}}\big)\mat{P}^{-1 \slash 2}\bigr) \\
    &\leq 1 + \norm{ \smash{\mat{P}^{-1 \slash 2}} }^2
    \norm{\smash{ \mat{A} - \mat{\hat{A}} }}
    = 1 + \norm{\smash{ \mat{P}^{-1} }}
    \norm{ \smash{\mat{A} - \mat{\hat{A}} }}.
\end{align*}
Finally, using the fact that $\lVert \mat{A} - \mat{\hat{A}} \rVert \leq \tr\big( \mat{A} - \mat{\hat{A}} \big)$ and $\lambda_{\rm min}(\mat{P}) \geq \mu$, we find
\begin{equation}
\label{eq:last}
    \lambda_{\rm max}\bigl(\mat{P}^{-1 \slash 2} (\mat{A} + \mu \Id)\mat{P}^{-1 \slash 2}\bigr) \leq 1 + \tr\big(\mat{A} - \mat{\hat{A}}\big) / \mu.
\end{equation}
Combining the bound \eqref{eq:conclude} for the minimum eigenvalue with the bound \eqref{eq:last} for the maximum eigenvalue verifies the condition number bound \eqref{eq:condition}.

If $\operatorname{rank}_\mu(\mat A)$ is the actual rank of $\mat{A}$, then \RPCholesky recovers $\mat{A}$ perfectly after $\operatorname{rank}_\mu(\mat A)$ steps.
Otherwise, we apply the main \RPCholesky error bound \cite[Thm.~5.1]{chen2022randomly} with rank $\rho = \text{rank}_{\mu}(\mat{A})$
and approximation accuracy
$1 + \mu/\sum_{i > \text{rank}_{\mu}(\mat{A})} \lambda_i(\mat{A}) \geq 2$, which guarantees
\begin{equation} 
\label{eq:bound_error}
    \expect \bigl[\tr\bigl(\mat{A} - \mat{\hat{A}}\bigr)\bigr]
    \leq \sum_{i > \text{rank}_{\mu}(\mat{A})} \lambda_i(\mat{A}) + \mu
\end{equation}
for any $r$ satisfying 
\begin{equation*}
    r \geq \textup{rank}_\mu(\mat{A}) \,\bigl(1 + \log\bigl(\tr \mat{A} / \mu\bigr)\bigr).
\end{equation*}
By the definition of the $\mu$-tail rank (\Cref{def:d_tail}), the right-hand side of \eqref{eq:bound_error} is at most $2\mu$.
Hence, combining \eqref{eq:bound_error} with the condition number bound \eqref{eq:condition} guarantees
\begin{equation*}
    \expect\bigl[ \kappa \bigl(\mat{P}^{-1 \slash 2} \bigl(\mat{A} + \mu \Id\bigr)\mat{P}^{-1 \slash 2}\bigr)\bigr] \leq 3.
\end{equation*}
By Markov's inequality, it follows that 
\begin{equation*}
    \kappa \bigl(\mat{P}^{-1 \slash 2} \bigl(\mat{A} + \mu \Id\bigr)\mat{P}^{-1 \slash 2}\bigr) \leq 3 \slash \delta
\end{equation*}
with failure probability at most $\delta$.
We apply the CG error bound \eqref{eq:rate} to obtain
\begin{equation*}
    \frac{\lVert \vec{\beta}^{(t)} - \vec{\beta} \rVert_{\mat{A} + \mu \Id}}
    {\lVert \vec{\beta} \rVert_{\mat{A} + \mu \Id}}
    \leq 2 \Biggl(\frac{\sqrt{3 \slash \delta} - 1}{\sqrt{3 \slash \delta} + 1}\Biggr)^t
    \leq 2 \e^{-2t \sqrt{\delta \slash 3}}
    \leq 2 \e^{-t\sqrt{\delta}}
\end{equation*}
for each $t \geq 0$, with failure probability at most $\delta$, which is equivalent to the stated convergence guarantee for \RPCholesky preconditioning.
\hfill $\qedsymbol$

\subsection{Proof of \texorpdfstring{\Cref{thm:a_priori_2}}{Theorem 1.3}}
\label{sec:a_priori_2}

Cohen \cite[Thm.~4.2]{Coh16} showed that a sparse sign embedding with parameters
$d = \order(k \log (k/\delta))$ and $\zeta = \order(\log (k/\delta))$
satisfies the following oblivious subspace embedding property:
for any $k$-dimensional subspace $\set{X} \subseteq \real^N$, there is probability at least $1 - 
\delta$ that
\begin{equation} \label{eq:distortion}
    (1/2) \norm{\vec{v}}^2 \le \norm{\mat{\Phi}\vec{v}}^2
    \le (3/2) \norm{\vec{v}}^2 \quad \text{for every $\vec{v} \in \set{X}$}.
\end{equation}
We set $\set{X} = \textup{range}(\mat{A}(:, \set{S}))$
and consider the matrices 
\begin{equation*}
\mat{P} = \mat{A}(\set{S},:) \mat{\Phi}^\ast \mat{\Phi} \mat{A}(:,\set{S}) + \mu \mat{A}(\set{S},\set{S}) 
\quad \text{and} \quad \mat{M} =  \mat{A}(\set{S},:)\mat{A}(:,\set{S}) + \mu \mat{A}(\set{S},\set{S}).
\end{equation*}
By the subspace embedding property \eqref{eq:distortion},
\begin{align*}
    \lambda_{\rm max}\mleft(\mat{P}^{-1/2}\mat{M}\mat{P}^{-1/2}\mright) 
    &= \max_{\vec{z}\ne \vec{0}} \frac{\vec{z}^*\mat{M}\vec{z}}{\vec{z}^*\mat{P}\vec{z}}
    \le \max_{\vec{z}\ne \vec{0}} 
    \;\max\Biggl\{1,
    \frac{\norm{\mat{A}(:,\set{S})\vec{z}}^2}{\norm{\mat{\Phi}\mat{A}(:,\set{S})\vec{z}}^2} \Biggr\}
    \le 2; \\
    \lambda_{\rm min}\mleft(\mat{P}^{-1/2}\mat{M}\mat{P}^{-1/2}\mright) 
    &= \min_{\vec{z}\ne \vec{0}}
    \frac{\vec{z}^*\mat{M}\vec{z}}{\vec{z}^*\mat{P}\vec{z}}
    \ge \min_{\vec{z}\ne \vec{0}} 
    \;\min\Biggl\{1, \frac{\norm{\mat{A}(:,\set{S})\vec{z}}^2}{\norm{\mat{\Phi}\mat{A}(:,\set{S})\vec{z}}^2}\Biggr\}
    \ge \frac{2}{3}.
\end{align*}
Therefore, $\kappa(\mat{P}^{-1/2}\mat{M}\mat{P}^{-1/2}) \leq 3$, and we apply the CG error bound \eqref{eq:rate} to obtain
\begin{equation*}
    \frac{\lVert \vec{\hat{\beta}}^{(t)} - \vec{\hat{\beta}} \rVert_{\mat{M}}}
    {\lVert \vec{\hat{\beta}} \rVert_{\mat{M}}}
    \leq 2 \Biggl(\frac{\sqrt{3} - 1}{\sqrt{3} + 1}\Biggr)^t
    \leq 2 \e^{-t}
\end{equation*}
for each $t \geq 0$, which is equivalent to the stated performance guarantee.\hfill \qedsymbol

\section{Conclusions}
\label{sec:conclusions}

We have developed two randomized preconditioners for kernel ridge regression problems with a moderate to large number of data points ($10^4\leq N\leq10^7$).
Across our experiments, \RPCholesky and \krill preconditioning were consistently competitive with or more robust than the tested alternatives \cite{COCF16,FTU21,gardner2018gpytorch,wang2019exact,RCR17,RCCR19,meanti2020kernel}.

If the number of data points $N$ is not prohibitively large (say, $N \leq 10^5$--$10^6$), we recommend solving the full-data KRR equations using conjugate gradient with \RPCholesky preconditioning.
The \RPCholesky approach enables us to solve KRR problems in just $\order(N^2)$ operations when the kernel matrix eigenvalues decay at a sufficiently fast polynomial rate and the chosen regularization is not too small.
For example, on the QM9 data set, \RPCholesky preconditioning solves the full-data KRR system with $N=10^5$ training points in $60$ CG iterations.

Eigenvalue decay is the main requirement limiting the performance of \RPCholesky preconditioning, and this limitation becomes apparent when the regularization $\mu$ is small (see the right panel of \Cref{fig:exact-performance}).
Future work is needed to discover effective preconditioners for KRR problems with slow eigenvalue decay.
Methods based on hierarchical decompositions \cite{COCF16,AFG+16,CAS17,YLRB17} and sparse inverse approximations \cite{schafer2021sparse,huan2025sparse,zhao2024adaptive} are promising candidates for this purpose, but they may perform poorly in high dimensions \cite{AFG+16}.

If the number of data points is so large that $\order(N^2)$ operations is too expensive,
we recommend restricting the KRR equations to $k \ll N$ data centers
and solving the restricted KRR equations using conjugate gradient with \krill preconditioning.
The \krill approach converges to the desired accuracy in just $\order((N + k^2)k \log k)$ operations for any kernel matrix and regularization parameter $\mu$.
For example, on the SUSY data set containing $5\times10^6$ points, \krill preconditioning solves the restricted KRR system formed from $N=4.5\times10^6$ training points in just four CG iterations, even with a parameter setting $\mu/N = 2\times 10^{-10}$ that would be challenging for other algorithms.

With appropriate parameter choices and neglecting rounding errors, \krill is guaranteed to control the condition number regardless of the regularization and the kernel matrix eigenvalue decay.
Nonetheless, there remain opportunities to improve on \krill.
First, the existing analysis of sparse random sign embeddings \cite{Coh16}
does not specify appropriate values of $\zeta$ and $d$ to use in practice.
More analysis of sparse sign embeddings is needed to close the theory--practice gap; recent papers \cite{chenakkod2025optimalsubspaceembeddingsresolving,tropp2025comparisontheoremsminimumeigenvalue} make some progress in this direction.
A second potential concern with \krill is numerical stability;
see the recent papers
\cite{MNTW24,EMN24}.
As a practical remedy, we have addressed numerical stability
by adding a small shift to the regularizer, but this approach has not been analyzed rigorously.
Finally, it may be possible to develop more efficient algorithms for restricted KRR problems when the number of data centers is $k \gg N^{1/2}$.
The ideal conjugate gradient algorithm would require $\order(kN)$ operations, but \krill is more expensive due to the cost of factorizing and applying the $k \times k$ preconditioner.

\section*{Declarations}

\begin{itemize}
\item \textbf{Conflict of interest.} The authors have no competing interests.
\item \textbf{Funding statement.} Ethan N. Epperly acknowledges support from 
DOE CSGF 
DE-SC0021110.
Mateo D\'iaz, Joel A. Tropp, and Robert J. Webber acknowledge support from 
ONR
BRC 
N00014-18-1-2363, from 
NSF FRG
1952777, and from the Caltech 
Carver Mead New Adventures Fund.
Zachary Frangella acknowledges support from
NSF IIS-1943131, the ONR
Young Investigator Program, and the Alfred P. Sloan Foundation.
\item \textbf{Author contributions.} Robert J. Webber was the primary writer; Mateo D\'iaz, Ethan N. Epperly, Zachary Frangella, and Robert J. Webber performed the numerical experiments; and all authors (Mateo D\'iaz, Ethan N. Epperly, Zachary Frangella, Joel A. Tropp, and Robert J. Webber) contributed to editing the paper.
%Author contribution
\item \textbf{Acknowledgments.} We thank Misha Belkin, Tyler Chen, Jackie Lok, Riley Murray, Madeleine Udell, and Jorge Garza-Vargas for helpful conversations.
\end{itemize}

\begin{appendices}

\section{Blocked randomly pivoted Cholesky}

To ensure reproducibility, we provide pseudocode for blocked \RPCholesky in \Cref{alg:rpcholesky}. 

 \begin{algorithm}[H]
\caption{{Blocked \RPCholesky for psd low-rank approximation} \label{alg:rpcholesky}}
\begin{algorithmic}[1]
\Require Positive semidefinite matrix $\mat{A} \in \real^{N \times N}$, approximation rank $k$, block size $B$.
\Ensure Factor matrix $\mat{F} \in \real^{N \times k}$, index set $\set{S}$.
\State Initialize $\mat{F} \leftarrow \mat{0}_{N\times k}$, $\set{S} \leftarrow \emptyset$, $\vec{d} \leftarrow \diag(\mat{A})$, $i \leftarrow 0$
\While {$i<k$}
\State Sample i.i.d. indices $s_1, \dots, s_{\min \{B, k-i\}} \sim \vec{d}/\Call{Sum}{\vec{d}}$
\State $\set{S}' \leftarrow \textsc{Unique}(s_1, \dots, s_{\min \{B, k-i\}})$, $ \set{S} \leftarrow \set{S} \cup \set{S}'$
\State $\mat{G} \leftarrow \mat{A}(:, \set{S}') - \mat{F}(:, 1:i) \mat{F}(\set{S}', 1:i)^*$ \label{line:rpc_G}
\State $\mat{R} \leftarrow \textsc{Cholesky}(\mat{G}(\set{S}', :))$ 
\Comment{$\mat{G}(\set{S}', :) = \mat{R}^*\mat{R}$} \label{line:rpc_R}
\State $\mat{F}\left(:, (i+1):(i+|\set{S}'|)\right) \leftarrow \mat{G}\mat{R}^{-1}$
\State $\vec{d} \leftarrow \vec{d} - \textsc{SquaredRowNorms}\left(\mat{F}\left(:, (i+1):(i+|\set{S}'|)\right)\right)$
\State $i \leftarrow i + |\set{S}'|$
\EndWhile
\end{algorithmic}
\end{algorithm}

\section{Data sets} \label{sec:data}
In our numerical experiments, we use data from LIBSVM \cite{chang2011libsvm}, OpenML \cite{vanschoren2014openml}, and UCI \cite{Dua:2019}.
We also use the QM9 data set \cite{RDRv14,RvBR12} which can be found online at \url{https://doi.org/10.6084/m9.figshare.c.978904.v5}.
See the GitHub repository \url{https://github.com/eepperly/Robust-randomized-preconditioning-for-kernel-ridge-regression} for a script to download the data sets.

         \begin{table}[t]
     \centering
     \caption{Data sets used in our experiments.}
     \label{tab:datasets}
    \begin{tabular}{llll} \toprule 
         \textbf{Data set}               & \textbf{Dimension} $d$ & \textbf{Total data-set size} $N$ & \textbf{Source}
          \\ \midrule
          \multicolumn{4}{c}{{Testbed}}
          \\ \midrule
          \texttt{ACSIncome}            & 11                     & 1\,331\,600              & OpenML\\ 
          \texttt{Airlines\_DepDelay\_1M} & 9                      & 800\,000                 & OpenML\\ 
          \texttt{cod\_rna}              & 8                      & 59\,535                  & LIBSVM\\ 
          \texttt{COMET\_MC\_SAMPLE}      & 4                      & 71\,712                  & LIBSVM\\ 
          \texttt{connect\_4}            & 126                    & 54\,045                  & LIBSVM\\ 
          \texttt{covtype\_binary}       & 54                     & 464\,809                 & LIBSVM\\ 
          \texttt{creditcard}           & 29                     & 227\,845                 & OpenML\\ 
          \texttt{diamonds}             & 9                      & 43\,152                  & OpenML\\ 
          \texttt{HIGGS}                & 28                     & 500\,000                 & LIBSVM\\ 
          \texttt{hls4ml\_lhc\_jets\_hlf}  & 16                     & 664\,000                 & OpenML\\ 
          \texttt{ijcnn1}               & 22                     & 49\,990                  & LIBSVM\\ 
          \texttt{jannis}               & 54                     & 46\,064                  & OpenML\\ 
          \texttt{Medical\_Appointment}  & 18                     & 48\,971                  & OpenML\\ 
          \texttt{MNIST}                & 784                    &60\,000                   & OpenML\\
          \texttt{sensit\_vehicle}       & 100                    & 78\,823                  & LIBSVM\\ 
          \texttt{sensorless}           & 48                     & 58\,509                  & LIBSVM\\ 
          \texttt{volkert}              & 180                    &46\,648                   & OpenML\\ 
          \texttt{w8a}                  & 300                    & 49\,749                  & LIBSVM\\ 
          \texttt{YearPredictionMSD}    & 90                     & 463\,715                 & LIBSVM\\ 
          \texttt{yolanda}              & 100                    &320\,000                  & OpenML\\ 
          \midrule
          \multicolumn{4}{c}{{Additional data sets}}             \\ \midrule
          \texttt{QM9}              & 435                    &133\,728                  & \cite{RDRv14,RvBR12} \\ 
          \texttt{SUSY}              & 18                    &5\,000\,000                 & UCI
         \\\bottomrule
     \end{tabular}
 \end{table}

\FloatBarrier
 
\end{appendices}

\bibliography{references}

@inproceedings{alaoui2015fast,
    author = {Alaoui, Ahmed and Mahoney, Michael W},
    booktitle = {Proceedings of the 28th International Conference on Neural Information Processing Systems},
    title = {Fast Randomized Kernel Ridge Regression with Statistical Guarantees},
    year = {2015},
    url = {https://dl.acm.org/doi/10.5555/2969239.2969326},
}

@article{AFG+16,
  title = {Fast Direct Methods for {G}aussian Processes},
  author = {Ambikasaran, Sivaram and {Foreman-Mackey}, Daniel and Greengard, Leslie and Hogg, David W. and O'Neil, Michael},
  year = {2016},
  journal = {IEEE Transactions on Pattern Analysis and Machine Intelligence},
  volume = {38},
  number = {2},
  pages = {252--265},
  doi = {10.1109/TPAMI.2015.2448083}
}

@article{avron2017faster,
    author = {Avron, Haim and Clarkson, Kenneth L. and Woodruff, David P.},
    title = {Faster Kernel Ridge Regression Using Sketching and Preconditioning},
    journal = {SIAM Journal on Matrix Analysis and Applications},
    volume = {38},
    number = {4},
    pages = {1116--1138},
    year = {2017},
    doi = {10.1137/16M1105396}
}

@article{AMT10,
  title = {Blendenpik: {{Supercharging LAPACK}}'s Least-Squares Solver},
  shorttitle = {Blendenpik},
  author = {Avron, Haim and Maymounkov, Petar and Toledo, Sivan},
  year = {2010},
  journal = {SIAM Journal on Scientific Computing},
  volume = {32},
  number = {3},
  pages = {1217--1236},
  doi = {10.1137/090767911},
}

@inproceedings{bach2013sharp,
  title = {Sharp analysis of low-rank kernel matrix approximations},
  author = {Bach, Francis},
  booktitle = {Proceedings of the 26th Annual Conference on Learning Theory},
  pages = {185--209},
  year = {2013},
  volume = {30},
  series = {Proceedings of Machine Learning Research},
  url = {https://proceedings.mlr.press/v30/Bach13.html}
}

@article{baldi2014searching,
  title = {Searching for exotic particles in high-energy physics with deep learning},
  author = {Baldi, Pierre and Sadowski, Peter and Whiteson, Daniel},
  journal = {Nature Communications},
  volume = {5},
  number = {4308},
  pages = {1--9},
  year = {2014},
  publisher = {Nature Publishing Group},
  doi = {10.1038/ncomms5308}
}

@article{caponnetto2006optimal,
  author = {Caponnetto, Andrea and De Vito, Ernesto},
  title = {Optimal Rates for the Regularized Least-Squares Algorithm},
  journal = {Foundations of Computational Mathematics},
  volume = {7},
  number = {3},
  pages = {331--368},
  year = {2007},
  doi = {10.1007/s10208-006-0196-8},
}

@article{chang2011libsvm,
    author = {Chang, Chih-Chung and Lin, Chih-Jen},
    title = {{LIBSVM}: {A} Library for Support Vector Machines},
    year = {2011},
    publisher = {Association for Computing Machinery},
    address = {New York, NY, USA},
    volume = {2},
    number = {3},
    doi = {10.1145/1961189.1961199},
    journal = {ACM Transactions on Intelligent Systems and Technology},
}

@article{CAS17,
  author = {Jie Chen and Haim Avron and Vikas Sindhwani},
  title = {Hierarchically Compositional Kernels for Scalable Nonparametric Learning},
  journal = {Journal of Machine Learning Research},
  year = {2017},
  volume = {18},
  number = {66},
  pages = {1--42},
  url = {https://jmlr.org/papers/v18/15-376.html}
}

@article{chen2022randomly,
  title = {Randomly Pivoted {{C}}holesky: {{Practical}} Approximation of a Kernel Matrix with Few Entry Evaluations},
  author = {Chen, Yifan and Epperly, Ethan N. and Tropp, Joel A. and Webber, Robert J.},
  year = {2025},
  journal = {Communications on Pure and Applied Mathematics},
  volume = {78},
  number = {5},
  pages = {995--1041},
  doi = {10.1002/cpa.22234},
}

@inproceedings{chenakkod2025optimalsubspaceembeddingsresolving,
  author = {Chenakkod, Shabarish and Derezi{\'n}ski, Micha{\l} and Dong, Xiaoyu},
  title = {Optimal Subspace Embeddings: Resolving {Nelson--Nguyen} Conjecture Up to Sub-Polylogarithmic Factors},
  booktitle = {Proceedings of the 2026 Annual ACM-SIAM Symposium on Discrete Algorithms (SODA)},
  pages = {4501--4559},
  year = {2026},
  doi = {10.1137/1.9781611978971.163},
}

@article{clarkson2017low,
    author = {Clarkson, Kenneth L. and Woodruff, David P.},
    title = {Low-Rank Approximation and Regression in Input Sparsity Time},
    year = {2017},
    volume = {63},
    number = {6},
    doi = {10.1145/3019134},
    journal = {Journal of the ACM},
}

@inproceedings{Coh16,
  author = {Michael B. Cohen},
  title = {Nearly Tight Oblivious Subspace Embeddings by Trace Inequalities},
  booktitle = {Proceedings of the Twenty-Seventh Annual ACM-SIAM Symposium on Discrete Algorithms},
  year = {2016},
  doi = {10.1137/1.9781611974331.ch21},
}

@inproceedings{COCF16,
  title = {Preconditioning Kernel Matrices},
  author = {Cutajar, Kurt and Osborne, Michael and Cunningham, John and Filippone, Maurizio},
  booktitle = {Proceedings of The 33rd International Conference on Machine Learning},
  year = {2016},
  url = {https://proceedings.mlr.press/v48/cutajar16.html},
}

@inproceedings{dai2014scalable,
    url = {https://dl.acm.org/doi/10.5555/2969033.2969166},
    author = {Dai, Bo and Xie, Bo and He, Niao and Liang, Yingyu and Raj, Anant and Balcan, Maria-Florina F and Song, Le},
    booktitle = {Proceedings of the 27th International Conference on Neural Information Processing Systems},
    title = {Scalable Kernel Methods via Doubly Stochastic Gradients},
    year = {2014}
}

@article{DBB21,
  author = {Deringer, Volker L. and Bartók, Albert P. and Bernstein, Noam and Wilkins, David M. and Ceriotti, Michele and Csányi, Gábor},
  title = {Gaussian Process Regression for Materials and Molecules},
  journal = {Chemical Reviews},
  volume = {121},
  number = {16},
  pages = {10073--10141},
  year = {2021},
  doi = {10.1021/acs.chemrev.1c00022},
}

@article{DM21,
  title = {Simpler Is Better: A Comparative Study of Randomized Pivoting Algorithms for {{CUR}} and Interpolative Decompositions},
  shorttitle = {Simpler Is Better},
  author = {Dong, Yijun and Martinsson, Per-Gunnar},
  year = {2023},
  month = aug,
  journal = {Advances in Computational Mathematics},
  volume = {49},
  number = {4},
  pages = {66},
  doi = {10.1007/s10444-023-10061-z},
}

@misc{Dua:2019,
  author = {Dua, Dheeru and Graff, Casey},
  title = {{UCI} Machine Learning Repository},
  year = {2019},
  institution = {University of California, Irvine, School of Information and Computer Sciences},
  url = {https://archive.ics.uci.edu/ml},
}

@article{EMN24,
  author = {Epperly, Ethan N. and Meier, Maike and Nakatsukasa, Yuji},
  title = {Fast Randomized Least-Squares Solvers Can Be Just as Accurate and Stable as Classical Direct Solvers},
  journal = {Communications on Pure and Applied Mathematics},
  volume = {79},
  number = {2},
  pages = {293--339},
  year = {2026},
  doi = {10.1002/cpa.70013},
}

@article{epperly2025embracerejectionkernelmatrix,
  author = {Epperly, Ethan N. and Tropp, Joel A. and Webber, Robert J.},
  title = {Embrace Rejection: Kernel Matrix Approximation by Accelerated Randomly Pivoted {Cholesky}},
  journal = {SIAM Journal on Matrix Analysis and Applications},
  volume = {46},
  number = {4},
  pages = {2527--2557},
  year = {2025},
  doi = {10.1137/24M1699048},
}

@article{FTU21,
  title = {Randomized {{Nystr{\"o}m}} Preconditioning},
  author = {Frangella, Zachary and Tropp, Joel A. and Udell, Madeleine},
  year = {2023},
  month = jun,
  journal = {SIAM Journal on Matrix Analysis and Applications},
  volume = {44},
  number = {2},
  pages = {718--752},
  doi = {10.1137/21M1466244},
}

@inproceedings{gardner2018gpytorch,
    author = {Gardner, Jacob and Pleiss, Geoff and Weinberger, Kilian Q and Bindel, David and Wilson, Andrew G},
    booktitle = {Proceedings of the 32nd International Conference on Neural Information Processing Systems},
    title = {{GPyTorch}: {Blackbox} Matrix-Matrix {Gaussian} Process Inference with {GPU} Acceleration},
    url = {https://dl.acm.org/doi/10.5555/3327757.3327857},
    year = {2018}
}

@article{huan2025sparse,
  author = {Huan, Stephen and Guinness, Joseph and Katzfuss, Matthias and Owhadi, Houman and Sch\"{a}fer, Florian},
  title = {Sparse Inverse {C}holesky Factorization of Dense Kernel Matrices by Greedy Conditional Selection},
  journal = {SIAM/ASA Journal on Uncertainty Quantification},
  volume = {13},
  number = {3},
  pages = {1649--1679},
  year = {2025},
  doi = {10.1137/23M1606253},
}

@misc{lacotte2021fast,
  author = {Lacotte, Jonathan and Pilanci, Mert},
  title = {Fast Convex Quadratic Optimization Solvers with Adaptive Sketching-Based Preconditioners},
  year = {2021},
  eprint = {2104.14101},
  archivePrefix = {arXiv},
  url = {https://arxiv.org/abs/2104.14101},
}

@article{LLR+22,
  title = {Learning new physics efficiently with nonparametric methods},
  doi = {10.1140/epjc/s10052-022-10830-y},
  author = {Letizia, Marco and Losapio, Gianvito and Rando, Marco and Grosso, Gaia and Wulzer, Andrea and Pierini, Maurizio and Zanetti, Marco and Rosasco, Lorenzo},
  journal = {The European Physical Journal C},
  volume = {82},
  number = {10},
  pages = {879},
  year = {2022},
  publisher = {Springer}
}

@inproceedings{ma2017diving,
    author = {Ma, Siyuan and Belkin, Mikhail},
    booktitle = {Proceedings of the 31st International Conference on Neural Information Processing Systems},
    title = {Diving into the shallows: {A} computational perspective on large-scale shallow learning},
    year = {2017},
    url = {https://dl.acm.org/doi/10.5555/3294996.3295135},
}

@article{MT20a,
  title = {Randomized Numerical Linear Algebra: {{Foundations}} and Algorithms},
  author = {Martinsson, Per-Gunnar and Tropp, Joel A.},
  year = {2020},
  journal = {Acta Numerica},
  volume = {29},
  pages = {403--572},
  doi = {10.1017/S0962492920000021},
}

@inproceedings{meanti2020kernel,
    author = {Meanti, Giacomo and Carratino, Luigi and Rosasco, Lorenzo and Rudi, Alessandro},
    booktitle = {Proceedings of the 34th International Conference on Neural Information Processing Systems},
    title = {Kernel Methods Through the Roof: {H}andling Billions of Points Efficiently},
    year = {2020},
    url = {https://dl.acm.org/doi/10.5555/3495724.3496932},
}

@misc{meier2022randomized,
  author = {Meier, Maike and Nakatsukasa, Yuji},
  title = {Randomized Algorithms for {Tikhonov} Regularization in Linear Least Squares},
  year = {2022},
  eprint = {2203.07329},
  archivePrefix = {arXiv},
  url = {https://arxiv.org/abs/2203.07329},
}

@article{MNTW24,
  title = {Are sketch-and-precondition least squares solvers numerically stable?},
  author = {Meier, Maike and Nakatsukasa, Yuji and Townsend, Alex and Webb, Marcus},
  journal = {SIAM Journal on Matrix Analysis and Applications},
  volume = {45},
  number = {2},
  pages = {905--929},
  year = {2024},
  publisher = {SIAM},
  doi = {10.1137/23M1551973},
}

@article{meng2014lsrn,
  author = {Meng, Xiangrui and Saunders, Michael A. and Mahoney, Michael W.},
  title = {{LSRN}: {A} Parallel Iterative Solver for Strongly Over- or Underdetermined Systems},
  journal = {SIAM Journal on Scientific Computing},
  volume = {36},
  number = {2},
  pages = {C95--C118},
  year = {2014},
  doi = {10.1137/120866580},
}

@techreport{murray2022,
  author = {Murray, Riley and Demmel, James and Mahoney, Michael W. and Erichson, N. Benjamin and Melnichenko, Maksim and Malik, Osman Asif and Grigori, Laura and Luszczek, Piotr and Derezi{\'n}ski, Micha{\l} and Lopes, Miles E. and Liang, Tianyu and Luo, Hengrui and Dongarra, Jack},
  title = {Randomized Numerical Linear Algebra: A Perspective on the Field with an Eye to Software},
  institution = {EECS Department, University of California, Berkeley},
  number = {UCB/EECS-2023-19},
  year = {2023},
  url = {https://www2.eecs.berkeley.edu/Pubs/TechRpts/2023/EECS-2023-19.html},
}

@inproceedings{musco2017recursive,
    title = {Recursive sampling for the {N}ystr\"om method},
    author = {Musco, Cameron and Musco, Christopher},
    url = {https://dl.acm.org/doi/10.5555/3294996.3295140},
    volume = {30},
    year = {2017},
    booktitle = {Proceedings of the 31st International Conference on Neural Information Processing Systems},
}

@inproceedings{ozaslan2019iterative,
    author = {Ozaslan, Ibrahim Kurban and Pilanci, Mert and Arikan, Orhan},
    booktitle = {IEEE International Conference on Acoustics, Speech and Signal Processing},
    title = {Iterative {H}essian Sketch with Momentum},
    year = {2019},
    doi = {10.1109/ICASSP.2019.8682720},
}

@article{pilanci2017newton,
    author = {Pilanci, Mert and Wainwright, Martin J.},
    title = {Newton Sketch: {A} Near Linear-Time Optimization Algorithm with Linear-Quadratic Convergence},
    journal = {SIAM Journal on Optimization},
    volume = {27},
    number = {1},
    pages = {205--245},
    year = {2017},
    doi = {10.1137/15M1021106},
}

@article{RDRv14,
  title = {Quantum Chemistry Structures and Properties of 134 Kilo Molecules},
  author = {Ramakrishnan, Raghunathan and Dral, Pavlo O. and Rupp, Matthias and {von Lilienfeld}, O. Anatole},
  year = {2014},
  journal = {Scientific Data},
  volume = {1},
  number = {1},
  pages = {140022},
  doi = {10.1038/sdata.2014.22}
}

@article{RT08,
  title = {A Fast Randomized Algorithm for Overdetermined Linear Least-Squares Regression},
  author = {Rokhlin, Vladimir and Tygert, Mark},
  year = {2008},
  journal = {Proceedings of the National Academy of Sciences},
  volume = {105},
  number = {36},
  pages = {13212--13217},
  publisher = {{Proceedings of the National Academy of Sciences}},
  doi = {10.1073/pnas.0804869105},
}

@article{RvBR12,
  title = {Enumeration of 166 Billion Organic Small Molecules in the Chemical Universe Database {GDB-17}},
  author = {Ruddigkeit, Lars and {van Deursen}, Ruud and Blum, Lorenz C. and Reymond, Jean-Louis},
  year = {2012},
  journal = {Journal of Chemical Information and Modeling},
  volume = {52},
  number = {11},
  pages = {2864--2875},
  doi = {10.1021/ci300415d},
}

@inproceedings{RCCR19,
    author = {Rudi, Alessandro and Calandriello, Daniele and Carratino, Luigi and Rosasco, Lorenzo},
    booktitle = {Proceedings of the 32nd International Conference on Neural Information Processing Systems},
    title = {On Fast Leverage Score Sampling and Optimal Learning},
    volume = {31},
    year = {2018},
    url = {https://dl.acm.org/doi/10.5555/3327345.3327470},
}

@inproceedings{RCR17,
    title = {{FALKON}: {A}n Optimal Large Scale Kernel Method},
    booktitle = {Proceedings of the 31st International Conference on Neural Information Processing Systems},
    author = {Rudi, Alessandro and Carratino, Luigi and Rosasco, Lorenzo},
    year = {2017},
    url = {https://dl.acm.org/doi/10.5555/3294996.3295145},
}

@book{Saa03,
    author = {Saad, Yousef},
    title = {Iterative Methods for Sparse Linear Systems},
    publisher = {Society for Industrial and Applied Mathematics},
    year = {2003},
    doi = {10.1137/1.9780898718003},
    edition   = {2},
    address   = {Philadelphia, PA},
}

@article{schafer2021sparse,
  author = {Sch\"{a}fer, Florian and Katzfuss, Matthias and Owhadi, Houman},
  title = {Sparse {C}holesky Factorization by {K}ullback--{L}eibler Minimization},
  journal = {SIAM Journal on Scientific Computing},
  volume = {43},
  number = {3},
  pages = {A2019-A2046},
  year = {2021},
  doi = {10.1137/20M1336254},
}

@inproceedings{smola2000sparse,
    title = {Sparse greedy {Gaussian} process regression},
    author = {Smola, Alex and Bartlett, Peter},
    booktitle = {Proceedings of the 13th International Conference on Neural Information Processing Systems},
    volume = {13},
    year = {2000},
    url = {https://dl.acm.org/doi/10.5555/3008751.3008838},
}

@article{STR+19,
  title = {Chemical Diversity in Molecular Orbital Energy Predictions with Kernel Ridge Regression},
  author = {Stuke, Annika and Todorovi{\'c}, Milica and Rupp, Matthias and Kunkel, Christian and Ghosh, Kunal and Himanen, Lauri and Rinke, Patrick},
  year = {2019},
  journal = {The Journal of Chemical Physics},
  volume = {150},
  number = {20},
  pages = {204121},
  publisher = {{American Institute of Physics}},
  doi = {10.1063/1.5086105}
}

@article{tropp2025comparisontheoremsminimumeigenvalue,
  author = {Tropp, Joel A.},
  title = {Comparison Theorems for the Minimum Eigenvalue of a Random Positive-Semidefinite Matrix},
  journal = {Communications of the American Mathematical Society},
  volume = {6},
  number = {16},
  pages = {844--895},
  year = {2026},
  doi = {10.1090/cams/72},
}

@article{vanschoren2014openml,
  author = {Vanschoren, Joaquin and {van Rijn}, Jan N. and Bischl, Bernd and Torgo, Luis},
  title = {{OpenML}: Networked Science in Machine Learning},
  journal = {ACM SIGKDD Explorations Newsletter},
  volume = {15},
  number = {2},
  pages = {49--60},
  year = {2014},
  doi = {10.1145/2641190.2641198},
}

@inproceedings{wang2019exact,
    title = {Exact {Gaussian} processes on a million data points},
    author = {Wang, Ke and Pleiss, Geoff and Gardner, Jacob and Tyree, Stephen and Weinberger, Kilian Q and Wilson, Andrew Gordon},
    booktitle = {Proceedings of the 33rd International Conference on Neural Information Processing Systems},
    year = {2019},
    url = {https://dl.acm.org/doi/10.5555/3454287.3455599},
}

@article{WM21,
  author = {Westermayr, Julia and Marquetand, Philipp},
  title = {Machine Learning for Electronically Excited States of Molecules},
  journal = {Chemical Reviews},
  volume = {121},
  number = {16},
  pages = {9873--9926},
  year = {2021},
  doi = {10.1021/acs.chemrev.0c00749}
}

@inproceedings{YLM+12,
    author = {Yang, Tianbao and Li, Yu-feng and Mahdavi, Mehrdad and Jin, Rong and Zhou, Zhi-Hua},
    booktitle = {Proceedings of the 25th International Conference on Neural Information Processing Systems},
    title = {Nystr\"{o}m Method vs Random {Fourier} Features: {A} Theoretical and Empirical Comparison},
    year = {2012},
    url = {https://dl.acm.org/doi/10.5555/2999134.2999188},
}

@inproceedings{YLRB17,
    author = {Yu, Chenhan D. and Levitt, James and Reiz, Severin and Biros, George},
    title = {Geometry-Oblivious {FMM} for Compressing Dense {SPD} Matrices},
    year = {2017},
    doi = {10.1145/3126908.3126921},
    booktitle = {Proceedings of the International Conference for High Performance Computing, Networking, Storage and Analysis},
}

@article{zhao2024adaptive,
  author = {Zhao, Shifan and Xu, Tianshi and Huang, Hua and Chow, Edmond and Xi, Yuanzhe},
  title = {An Adaptive Factorized {N}ystr\"om Preconditioner for Regularized Kernel Matrices},
  journal = {SIAM Journal on Scientific Computing},
  volume = {46},
  number = {4},
  pages = {A2351-A2376},
  year = {2024},
  doi = {10.1137/23M1565139},
}

@inproceedings{huang2024quasi,
  author    = {Huang, Zhen and Sun, Jiajin and Huang, Yian},
  title     = {Quasi-{M}onte {C}arlo Features for Kernel Approximation},
  booktitle = {Proceedings of the 41st International Conference on
               Machine Learning},
  series    = {Proceedings of Machine Learning Research},
  volume    = {235},
  pages     = {20134--20165},
  year      = {2024},
  publisher = {PMLR},
  address = {Vienna, Austria},
  url       = {https://proceedings.mlr.press/v235/huang24w.html}
}

@inproceedings{shi2023kernelanchor,
  author    = {Shi, Wenqi and Xu, Wenkai},
  title     = {Learning Nonlinear Causal Effect via Kernel Anchor Regression},
  booktitle = {Proceedings of the Thirty-Ninth Conference on
               Uncertainty in Artificial Intelligence},
  series    = {Proceedings of Machine Learning Research},
  volume    = {216},
  pages     = {1942--1952},
  year      = {2023},
  publisher = {PMLR},
  address = {Pittsburgh, PA, USA},
  url       = {https://proceedings.mlr.press/v216/shi23a.html}
}

\end{document}